 \documentclass[12pt]{amsart}
\usepackage{amsmath,amssymb}
\usepackage{amsfonts}
\usepackage{amsthm}
\usepackage{latexsym}
\usepackage{graphicx}

\def\e{\epsilon}
\def\lf{\left}
\def\ri{\right}

\def\wt{\widetilde}

\def\p{\partial}

\def\la{\langle}
\def\ra{\rangle}

\def\K{K\"ahler  }
\def\KR{K\"ahler-Ricci }
\def\KRF{K\"ahler-Ricci flow }

\def\be{\begin{equation}}
\def\ee{\end{equation}}

\def\lf{\left}
\def\ri{\right}

\def\e{\epsilon}

\def\wt{\widetilde}

\def\p{\partial}

\def\wh{\widehat}

\def\wt{\widetilde}

\def\p{\partial}

\def\p{\partial}

\def\KRF{K\"ahler-Ricci  flow }

\def\D{\displaystyle}

\def\supp{\mbox{supp}}

\def\cl{\overline}

\def\la{\langle}
\def\ra{\rangle}
\def\p{\partial}
\def\K{K\"ahler}

\newtheorem{thm}{Theorem}[section]
\newtheorem{ass}{Assumption}
\newtheorem{lem}{Lemma}[section]
\newtheorem{prop}{Proposition}[section]
\newtheorem{cor}{Corollary}[section]
\theoremstyle{definition}
\newtheorem{defn}{Definition}[section]
\theoremstyle{remark}

\newtheorem{rem}{Remark}[section]
\numberwithin{equation}{section}

\begin{document}
\title{Pseudolocality for the Ricci flow and applications}

\author{Albert Chau$^1$}
\thanks{$^1$Research
partially supported by NSERC grant no. \# 327637-06}

\address{Waterloo University, Department of Pure Mathematics,
  200 University avenue, Waterloo, ON N2L 3G1, CANADA}
\email{a3chau@math.uwaterloo.ca}

\author{Luen-Fai Tam$^2$}

\thanks{$^2$Research
partially supported by Earmarked Grant of Hong Kong \#CUHK403005}

\address{Department of Mathematics, The Chinese University of Hong Kong,
Shatin, Hong Kong, China.} \email{lftam@math.cuhk.edu.hk}
\author{ Chengjie Yu}\address{Department of Mathematics, The Chinese University of Hong Kong,
Shatin, Hong Kong, China.} \email{cjyu\@math.cuhk.edu.hk}

\renewcommand{\subjclassname}{%
  \textup{2000} Mathematics Subject Classification}
\subjclass[2000]{Primary 53C44; Secondary 58J37, 35B35}

\date{December,  2006}



\begin{abstract}
In \cite{P1}, Perelman established a differential Li-Yau-Hamilton
(LYH) type inequality for fundamental solutions of the conjugate
heat equation corresponding to the Ricci flow on compact manifolds
(also see \cite{N2}).  As an application of the LYH inequality,
Perelman proved a pseudolocality result for the Ricci flow on
compact manifolds. In this article we provide the details for the
proofs of these results in the case of a complete non-compact
Riemannian manifold. Using these results we prove that under
certain conditions, a finite time singularity of the Ricci flow
must form within a compact set.  We also prove a long time
existence result for the \KRF flow on complete non-negatively
curved \K\  manifolds.
\end{abstract}

\maketitle \markboth{Albert Chau, Luen-Fai Tam and Chengjie Yu}
{Pseudolocality and applications}

\section{Introduction}

 In this article we consider the Ricci flow

\begin{equation}\label{ie1}
\frac{\p}{\p t}g_{ij}=-2R_{ij}
\end{equation}
on a complete non-compact Riemannian manifold $(M, g)$, and the
heat equation and  conjugate heat equation
\begin{equation}\label{ie3}
\frac{\partial u}{\partial t}-\Delta^t u=0
\end{equation}

 \begin{equation}\label{ie2}
\frac{\partial u}{\partial t}+\Delta^t u-Ru=0
\end{equation}
where $\Delta^t$ denotes the Laplacian operator with respect to a
solution $g(t)$ to (\ref{ie1}), and $R(t)$ is the scalar curvature
of $g(t)$.  Notice that if $g(t)$ is defined on an interval $[0,
T]$ and we let $\tau=T-t$, then (\ref{ie2}) defines a strictly
parabolic equation on $M$ with respect to $\tau \in [0, T]$.

The conjugate heat equation corresponding to the Ricci flow was
considered  in \cite{P1}, and there Perelman established a
differential Li-Yau-Hamilton (LYH) type inequality for its
fundamental solutions (\cite{P1}; Corollary 9.3) on compact
manifolds. The proof was sketched in \cite{P1} and a detailed
proof was given by Ni in \cite{N2}.  As an application of the LYH
inequality, Perelman proved a pseudolocality result for the Ricci
 flow on compact manifolds (\cite{P1}; Theorem 10.1), which
  basically states that regions of large amounts of curvature
  cannot instantly
 affect almost Euclidean regions under the Ricci flow. For more
 details of the proof, see \cite{CLN,KL,ST}.

 In this article we verify these results, including the LHY hamilton inequality and pseudolocality,
 in the case of complete non-compact Riemannian manifolds.
   We basically follow the original steps described in \cite{P1}
    as well as those in \cite{CLN,KL,ST,N2}.

 Our motivation to generalize Perelman's results mentioned above is to
 study long time existence of Ricci flow and \KR flow on complete
 noncompact manifold. Using the result on pseudolocality, we obtain the following:

\begin{thm}\label{it1}
Let $(M^n, g)$ be a complete noncompact Riemannian manifold with
injectivity radius bounded away from zero such that
$$|Rm|(x)\to0$$ as $x\to \infty$. Let    $(M, g(t))$ be the
corresponding $maximal$ solution to the Ricci flow (\ref{lte1}) on
$M\times[0, T)$.  Then either $T=\infty$ or there exists some
compact $S \subset M$ with the property that $|Rm(x, t)|$ is
bounded on $(M \setminus S) \times [0, T)$.
\end{thm}

The conditions are satisfied if $M$ is an asymptotically flat
manifold for example.  As a Corollary to Theorem \ref{it1} we also have

\begin{cor}
Suppose $T < \infty$ in Theorem \ref{it1}.  Then $Rm(x, T) \to 0$ as $x\to \infty$ in the sense that: given any $\epsilon>0$, we may choose $S$ such that $|Rm(x, t)| \leq \epsilon$ for all $(x, t)\in S^c \times [0, T)$.
\end{cor}

Combining Theorem \ref{it1} with the results in \cite{NT2}, we have the
following result on the long time existence of \KR on complete
noncompact K\"ahler manifolds with nonnegative holomorphic
bisectional curvature.

\begin{thm}\label{ic1}
Let $(M^n, g_0)$ be a complete non-compact \K  \ manifold with
non-negative holomorphic bisectional curvature with injectivity
radius bounded away from zero such that $$|Rm|(x)\to0$$ as $x\to
\infty$.     Then   the \KR flow
\begin{equation}\label{KRflow}
\frac{\p}{\p t}g_{i\bar j}=-R_{i\bar j}
\end{equation}
with initial data $g_0$ has a long time solution $g(t)$ on
$M\times [0, \infty)$.
\end{thm}
The \KR flow is an important tool to study uniformization of
complete noncompact K\"ahler manifolds with nonnegative
holomorphic bisectional curvature, see \cite{Sh2,CT3,CT4} for
example. In \cite{Sh2} (see also \cite{NT2}), it was proved that
if $(M^n, g_0)$ is a complete non-compact \K manifold with
non-negative and bounded holomorphic bisectional curvature, and if
the scalar curvature satisfies:
\begin{equation}\label{scalarcurvaturedecay1}
\frac1{V_x(r)}\int_{B_x(r)}R\le \frac{C}{1+r^\theta}
\end{equation}
for some $C, \theta>0$ for all $x$ and $r$, then (\ref{KRflow})
has long time solution. By the result in \cite{NT1},
(\ref{scalarcurvaturedecay1}) is true for $\theta=1$, at least for
simply connected $M$ and where the constant $C$ which may  depend on $x$.
It is unclear whether (\ref{scalarcurvaturedecay1}) is true in
general with $C$ being independent of $x$ except for the case of
maximal volume growth, see \cite{N4}.

In order to prove the LYH type differential inequality for the fundamental solution  of
(\ref{ie3}),  we need to obtain estimates for the fundamental solution together with some
gradient estimates for positive solutions of (\ref{ie3}) and (\ref{ie2}). In case the manifold
is compact, results have been obtained by Zhang, Kuang-Zhang \cite{Zhang,Kuang-Zhang} and Ni
\cite{N2}. Some estimates are also obtained for complete manifolds with nonnegative Ricci
curvature by Ni \cite{N3}. We consider the case that the manifold is complete, non-compact, and
has bounded curvature. The results may have independent interest.

 The paper is organized as follows.  In every section,  our results are
 obtained on a complete non-compact Riemannian manifold.
 In \S 2 -\S 4 we establish  some basic estimates for positive solutions
 of the conjugate heat equation associated to a general evolution
 (\ref{evolution1}) of a Riemannian metric.  In \S 5 and \S 6 we establish
 estimates for fundamental solutions of this conjugate heat equation.
 In \S 7 we apply our previous estimates to establish the LYH inequality
  for the fundamental solution of the conjugate heat equation associated to
  the Ricci flow (\ref{ie1}).  Our steps in this section basically follow
  the steps in \cite{N2}
   \footnote{Our proof does not use the reduced distance $L(y, \tau)$
   associated to the Ricci flow, introduced in \cite{P1})}.
   In \S 8 we establish pseudolocality for the Ricci flow (\ref{ie1})
   on complete non-compact Riemannian manifolds.
    In particular, we show that Theorem 10.1 in \cite{P1}
    holds in the non-compact case.
    In \S 9 we prove Theorem \ref{it1} and Theorem \ref{ic1}.

\section{an integral estimate}
In this section, we will modify the arguments by Grigor'yan
\cite{Grigor'yan-1} a little bit to obtain an integral estimate
for solution of (\ref{eqn-schordinger}). The proof is basically
the same as in \cite{Grigor'yan-1}.

 Let $\{g(t)|t\in[0,T]\}$ be a smooth family of complete
Riemannian metrics on $M^n$ such that $g(t)$ satisfies:
\begin{equation}\label{evolution1}
\frac{\p}{\p t}g_{ij}(x,t)=2h_{ij}(x,t)
\end{equation}
on  $M\times [0,T]$, where $h_{ij}(x,t)$ is a smooth family of
symmetric tensors.

Consider the equation:
\begin{equation}\label{eqn-schordinger}
\frac{\partial u}{\partial t}-\Delta^t u+qu=0
\end{equation}
where $\Delta^t$ denotes the Laplacian operator with respect to
$g(t)$ and $q$ is a smooth function on $M\times [0,T]$

Let us make the following assumptions:

\begin{itemize}
    \item  [\textbf{ (A1)} ] $||h||, ||\nabla^th||$
     are uniformly bound on space-time, where the norm is taken with respect to $t$.
     \\
\item  [\textbf{ (A2)}] The sectional curvatures of the metrics
$g(t)$ are uniformly bounded on space-time.\\
\item  [\textbf{ (A3)}]$|q|,||\nabla^t q||,|\Delta^t q|$ are
uniformly bounded on space-time.
\end{itemize}

Let $H(t)$ be the trace of $h_{ij}(t)$ with respect to $g(t)$.

\begin{defn}\label{defn-regular}
  Let $f$ be a positive function on $(0,T]$.  $f$ is said to be
  regular with the constants $\gamma>1$ and $A\geq 1$, if

  (i) $f$ is increasing, and

  (ii) $\D\frac{f(s)}{f(s/\gamma)}\leq A\frac{f(t)}{f(t/\gamma)}$ for
  any $0<s\leq t\leq T$.
\end{defn}
\begin{lem}\label{lem-integral-estimate}
Let  $\Omega$ be a relative compact domain  of $M$ with smooth
boundary and let $K$ be a compact set with $K\subset\subset
\Omega$. Let $u$ be any solution to the problem:
\begin{equation}\label{eqn-dirichlet}
\left\{\begin{array}{l}u_t-\Delta^t u+qu=0,\  \text{in
$\Omega\times[0,T]$}\\u\big|_{\partial
\Omega\times[0,T]}=0\\\supp\, u(\cdot,0)\subset
K.\end{array}\right.
\end{equation}
 Let $f$
be a regular function with the constants $\gamma$ and $A$. Suppose
\begin{equation*}
\int_{\Omega}u^2dV_t\leq\frac{1}{f(t)}
\end{equation*}
for any $t>0$. Then there is a positive constant $C$ depending
only on $\gamma$, the uniform  upper bound of $|q|$ and $|H|$,
and a positive constant $D$ depending only on  $T$, $\gamma$ and
the uniform  upper bound of $||h||$, such that
\begin{equation*}
\int_{\Omega}u^2(x,t)e^{\frac{r^2(x,K)}{Dt}}dV_t\leq\frac{4A}{f(t/\gamma)}e^{Ct}
\end{equation*}
for any $t>0$, where $r(x,K)$ denotes the distance between $x$ and
$K$ with respect to the initial metric.
\end{lem}
\begin{proof}  The proof is almost the same as the proof of Theorem
2.1 in Grigoryan \cite{Grigor'yan-1}.

Let $C_2>0$ be a constant such that $|q|+1/2|H|\leq C_2$, and let
$v=e^{-C_2t}u$. Then $v$ satisfies
\begin{equation}\label{e1}
v_t-\Delta^t v+(C_2+q)v=0,
\end{equation}
and
\begin{equation*}
\int_{\Omega}v^2dV_t=e^{-2C_2t}\int_{\Omega}u^2dV_t\leq\frac{1}{f(t)e^{2C_2t}}:=\frac{1}{\tilde
f(t)}
\end{equation*}
where $\tilde f(t)=f(t)e^{2C_2t}$ is regular with constants
$\gamma$ and $A$.

 Now for any $R>0$, define
\begin{equation*}
d(x)=\left\{\begin{array}{ll}R-r(x,K)&x\in K^R\\0&x\not\in
K^R\end{array}\right.
\end{equation*}
where $K^R$ means the $R$-neighborhood of $K$ with respect to the
initial metric. Then $|\nabla^t d|\leq C_1$ uniformly on
space-time where $C_1$ depends on $T$ and the  upper bound of
$||h||$.  Then if we let
$\xi(x,s-t)=\frac{d^2(x)}{2C_{1}^2(t-s)}$ for $s>T$ fixed and $0<t\leq T<s$, we have
\begin{equation}\label{xi}
\frac{\p}{\p t}\xi+\frac{1}{2}
|\nabla^t\xi|^2=-\frac{d^2}{2C_{1}^2
(t-s)^2}+\frac{1}{2}\frac{d^2||\nabla^t d||^2}{C_1^{4}(t-s)^2}\leq
0,
\end{equation}
which combines with (\ref{e1}) to give
\begin{eqnarray*}
 & &\frac{d}{dt}\int_{\Omega}v^2e^{\xi}dV_t\\
 &=&\int_{\Omega}(2vv_te^{\xi}+v^2e^\xi\xi_t)dV_t+\int_{\Omega}v^2e^\xi HdV_t\\
 &\leq&2\int_{\Omega}v(\Delta^tv-(C_2+q)v)e^\xi
 dV_t-\frac{1}{2}\int_{\Omega}v^2e^\xi|\nabla^t\xi|^2dV_t-\int_{\Omega}Hv^2e^\xi
 dV_t\\
 &\leq&-2\int_{\Omega}|\nabla^tv|^2e^\xi dV_t-2\int_{\Omega}ve^\xi\la\nabla^tv,\nabla^t
 \xi\ra dV_t-\frac{1}{2}\int_{\Omega}v^2e^\xi|\nabla^t\xi|^2dV_t\\
 &\leq& 0,
 \end{eqnarray*}
 where we have used the fact that $v=0$ on $\p\Omega$,
 $2C_2+2q+H\ge0$ and (\ref{xi}).

 It now follows from STEP1 and STEP2 of the proof of Theorem 2.1 in Griyor'yan
\cite{Grigor'yan-1}, that there exists a positive constant $D>0$ depending
on $C_1$ and $\gamma$, such that
\begin{equation*}
 \int_{\Omega}v^2e^{\frac{r^2(x,K)}{Dt}}dV_t\leq\frac{4A}{\tilde
 f(t/\gamma)}.
 \end{equation*}
Thus,
\begin{equation*}
\int_{\Omega}u^2e^{\frac{r^2(x,K)}{Dt}}dV_t\leq\frac{4A}{f(t/\gamma)}
e^{2C_2(1-\frac{1}{\gamma})t}.\quad
\end{equation*}
\end{proof}

\section{A mean value inequality}

In this section, we prove the following lemma which will be used
to estimate the fundamental solution  of (\ref{backward}).
\begin{lem}\label{lm-mean-value-inequality}
Let $u$ be a positive sub-solution of equation
(\ref{eqn-schordinger}) on $\Omega\times [0,T]$ where $\Omega$ is
a domain in $M$. Moreover, suppose that there is a complete
Riemannian metric $\tilde g$ on $M$ with Ricci curvature bounded
from below by $-k$ with $k\geq 0$, such that
\begin{equation*}
\frac{1}{C_{0}}\tilde g\leq g(0)\leq C_{0}\tilde g
\end{equation*}
in $\Omega$ for some $C_0>0$. Let $\tilde Q_r(x,t):=\tilde
B_x(r)\times (t-r^2,t]$ whenever it is well defined, where $\tilde
B_x(r)$ means the ball of radius $r$ with respect to $\tilde g$.
Then for any $(x,t)\in\Omega\times(0,T]$ and $r>0$ such that
$\tilde Q_{2r}(x,t)\subset\subset\Omega\times[0,T]$,
\begin{equation*}
\sup_{\tilde Q_r(x,t)}u\leq\frac{Ce^{At+B\sqrt{k}r}}{r^2\tilde
V_x(r)}\int_{\tilde Q_{2r}(x,t)}ud\tilde Vds
\end{equation*}
where $A$ depends only on the the   upper bounds of $|q|$ and
$|H|$ on $\Omega$, $B$ depends only on $n$, and $C$ depends only
on  $C_0$, $n$, $T$ and the uniform  upper bound  of $|h|$ on
$\Omega$. The notations $\tilde V_x(r)$ and $d\tilde V$ denote the
volumes with respect to $\tilde g$.
\end{lem}
\begin{proof} The proof is almost the same as in Zhang [\cite{Zhang}; $\S$ 5]. We
just emphasis on our modifications. For details, please refer to
the paper.

Let $\sigma$ be in $(1,2]$. Let $\phi$ be a smooth function on
$[0,\infty)$, such that: i) $\phi=1$ on $[0,r]$, ii) $\phi=0$ on
$[\sigma r, \infty)$ and iii) $-\frac{2}{(\sigma-1)r}\leq\phi'\leq
0.$

Let $\eta$ be a smooth function on $[0,\infty)$, such that: i)
$\eta=0$ on $[0,t-\sigma^2 r^2]$, ii) $\eta=1$ on $[t-r^2,\infty)$
and iii) $0\leq \eta'\leq\frac{2}{(\sigma-1)^2r^2}$.

Let $\psi(y,s)=\phi(\tilde r(x,y))\eta(s)$ where $x$ is fixed and
$\tilde r(x,y)$ means the distance function of $\tilde g$. Then
$\supp \psi(\cdot,s)\subset \cl{\tilde B_x(2r)}\subset \Omega$.

Let $v=e^{-C_1t}u$, where $C_1$ is some positive constant to be
determined. Then
\begin{equation*}
v_t-\Delta^tv+(C_1+q)v\leq 0.
\end{equation*}

For any $p\geq 1$,
\begin{equation}\label{meanvalue1}
\frac{\partial v^p}{\partial t}-\Delta^t(v^p)+p(C_1+q)v^p\leq 0.
\end{equation}

Let $w=v^p$. Let $t'$ be any real number in $[t-r^2,t]$. Multiply
$w\psi^2$ to the inequality above and integrate. We get
\begin{equation}\label{meanvalue2}
-\int_{t-\sigma^2r^2}^{t'}\int_{M}\psi^2w\Delta^swdV_sds+p\int_{t-\sigma^2r^2}^{t'}\int_M(C_1+q)\psi^2w^2dV_sds
\end{equation}
\begin{equation}\nonumber
\leq -\int_{t-\sigma^2r^2}^{t'}\int_M\psi^2ww_sdV_sds\\.
\end{equation}

Integrating by parts in the first term on the left in the above
inequality, we get
\begin{equation}\label{meanvalue3}
\begin{split}
 -\int_{t-\sigma^2r^2}^{t'}\int_{M}\psi^2w\Delta^swdV_sds&
 = \int_{t-\sigma^2r^2}^{t'}\int_{M}\langle\nabla^s(\psi^2w),\nabla^sw\rangle_sdV_sds\\
&=\int_{t-\sigma^2r^2}^{t'}\int_{M}(|\nabla^s(\psi
w)|^2-|\nabla^s\psi|^2w^2)dV_sds.
\end{split}
\end{equation}

Moreover, the second term on the left in (\ref{meanvalue2}) is:
\begin{equation}\label{meanvalue4}
\begin{split}-\int_{t-\sigma^2r^2}^{t'}&\int_M\psi^2ww_sdV_sds\\
&= -\int_{t-\sigma^2r^2}^{t'}\int_M\psi^2ww_se^FdVds\\
&=-\int_{M}\int_{t-\sigma^2r^2}^{t'}\frac12\lf[(\psi^2w^2e^F)_s-
(\psi^2e^F)_sw^2\ri]dsdV\\
&=-\frac{1}{2}\int_{M}w^2\psi^2dV_{t'}+\int_{t-\sigma^2r^2}^{t'}\int_M\psi\psi_sw^2dV_sds\\
&\hspace{12pt}+\frac{1}{2}\int_{t-\sigma^2r^2}^{t'}\int_M
H\psi^2w^2dV_sds \end{split}
\end{equation}
where $F$ is such that $e^FdV=dV_s$ and $dV$ is the volume element
for $g(0)$.

 Choose $C_1$ large enough depending only on the
uniform  upper bounds of $|q|$ and $|H|$ on $\Omega$. Then by
(\ref{meanvalue2})--(\ref{meanvalue4}), we have:
\begin{eqnarray*}
\begin{split}
\int_{t-\sigma^2r^2}^{t'}&\int_{M}|\nabla^s(\psi
w)|^2dV_sds+\frac{1}{2}\int_{M}w^2\psi^2dV_{t'}\\
&\leq
\int_{t-\sigma^2r^2}^{t'}\int_M(\psi\psi_s+|\nabla^s\psi|^2)w^2dV_sds.
\end{split}
\end{eqnarray*}

Note that
\begin{equation*}
\psi_s=\phi\eta'\leq\frac{2}{(\sigma-1)^2r^2}
\end{equation*}
Moreover
\begin{eqnarray*}
|\nabla^s\psi|^2=\eta^2(\phi')^2|\nabla^s(\tilde
r(x,\cdot))|^2\leq\frac{C_2}{(\sigma-1)^2r^2}
\end{eqnarray*}
where $C_2$ depends on $C_0$, $T$ and the uniformly upper bound of
$|h|$ on $\Omega$.

Hence
\begin{eqnarray*}
& &\int_{t-\sigma^2r^2}^{t'}\int_{\tilde B_x(\sigma
r)}|\nabla^s(\psi
w)|^2dV_sds+\frac{1}{2}\int_{M}w^2\psi^2dV_{t'}\\
&\leq&
\frac{C_3}{(\sigma-1)^2r^2}\int_{t-\sigma^2r^2}^{t'}\int_{\tilde
B_x(\sigma r)}w^2dV_sds
\end{eqnarray*}
where $C_3=\max\{2,C_2\}$. So,
\begin{equation*}
\int_{t-\sigma^2r^2}^{t}\int_{\tilde B_x(\sigma r)}|\tilde
\nabla(\psi w)|^2d\tilde Vds\leq
\frac{C_4}{(\sigma-1)^2r^2}\int_{t-\sigma^2r^2}^{t}\int_{\tilde
B_x(\sigma r)}w^2d\tilde Vds\ \ \mbox{and}
\end{equation*}
\begin{equation*}
\max_{t-r^2\leq t'\leq t}\int_{\tilde B_x(\sigma
r)}w^2\psi^2d\tilde V\leq
\frac{C_4}{(\sigma-1)^2r^2}\int_{t-\sigma^2r^2}^{t}\int_{\tilde
B_x(\sigma r)}w^2d\tilde Vds
\end{equation*}
where $C_4$ depends only on   $C_0$, $T$ and the uniform  upper
bound of $|h|$.

Now we can proceed as in the proof of Theorem 5.1 in \cite{Zhang}
with respect to the metric $\tilde g$. Applying the the Sobolev
inequality  in \cite{SC} with respect to $\tilde g$ and the Moser
iteration as in Zhang \cite{Zhang}, we get

\begin{equation*}
\sup_{\tilde Q_r(x,t)}v^2\leq\frac{C_{5}e^{C_6\sqrt k
r}}{r^2\tilde
V_x(r)}\Big(\frac{1}{\log\gamma}\Big)^{\frac{n+2}{2}}\int_{\tilde
Q_{\gamma r}(x,t)}v^2d\tilde Vds
\end{equation*}
for any $\gamma\in(1,2]$, where $C_5$ depends only on $C_0$, $T$
the uniform  upper bound of $|h|$ and $n$, and $C_6$ depends only
$n$.

By the trick on iteration of Li-Schoen \cite{Li-Schoen} as
mentioned in \cite{Zhang}, we get
\begin{equation*}
\sup_{\tilde Q_r(x,t)}v\leq\frac{C_{7}e^{C_8\sqrt kr}}{r^2\tilde
V_x(r)}\times\frac{1}{\log^{\frac{n+2}{2}}\gamma}\int_{\tilde
Q_{\gamma r}(x,t)}vd\tilde Vds
\end{equation*}
for any $\gamma\in (1,2]$, where $C_7$ depends only on$C_0$, $T$ the upper bound
of $|h|$ and $n$, and $C_8$ depends only on $n$.

In particular, let $\gamma=2$, we get
\begin{equation*}
\sup_{\tilde Q_r(x,t)}v\leq\frac{C_{9}e^{C_{8}\sqrt
{k}r}}{r^2\tilde V_x(r)}\int_{\tilde Q_{2r}(x,t)}vd\tilde Vds
\end{equation*}
where $C_9$ depends only on $C_0$, $T$ the uniform  upper bound of
$|h|$ and $n$. So
\begin{eqnarray*}
\sup_{\tilde Q_r(x,t)}u&\leq& e^{C_1t}\sup_{\tilde Q_r(x,t)} v\\
&\leq&\frac{C_{9}e^{C_1t+{C_{8}\sqrt kr}}}{r^2\tilde
V_x(r)}\int_{\tilde Q_{2r}(x,t)}ud\tilde V ds.
\end{eqnarray*}
\end{proof}
\section{A Li-Yau type gradient estimate}
In this section, we derive a Li-Yau type graident estimate which
will also be used in the estimates of the fundamental solution of
\ref{conjugateheat}.  We  basically follow the proof in
\cite{LY}. Let $g(t)$ be as in \S2.

\begin{lem}\label{lm-harnack}
Let $u$ be a positive solution to equation
(\ref{eqn-schordinger}). Then for any $\alpha>1$ and $\epsilon>0$,
there is a constant $C>0$ depending on $\alpha,\epsilon, n, T$,
the uniform   upper bounds of $|h|,|\nabla h|$,$|\nabla q|,|\Delta
q|,|Rc^t|$ and the bound of the sectional curvature at $t=0$, such
that
\begin{equation*}
\frac{|\nabla u|^2}{u^2}-\alpha\frac{u_t}{u}-\alpha q\leq
C+\frac{(n+\epsilon)\alpha^2}{2t}.
\end{equation*}
\end{lem}
\begin{proof} In the following $\nabla$ and $\Delta$ are understood to be time
dependent. For any smooth function $f$ on $M\times [0,T)$, at a
point with normal coordinates with respect to the metric $g(t)$,
we have
\begin{eqnarray*}
(\Delta f)_t&=&\Delta f_t-2h_{ij}f_{ij}-2h_{ik;i}f_k+H_if_i\\
(|\nabla f|^2)_t&=&2\la\nabla f_t,\nabla f\ra-2h(\nabla f,\nabla
f).
\end{eqnarray*}
Repeated indices mean summation.

 Let $f=\log u$. Then
\begin{equation}\label{eqn-f}
\Delta f-f_t=q-|\nabla f|^2.
\end{equation}

For $\alpha>1$ and $\e>0$, let $F=t(|\nabla f|^2-\alpha f_t-\alpha
R)$. Then in    normal coordinates
\begin{equation}\label{LaplacianF}
   \begin{split}
 \Delta F&=t\bigg[2\sum_{ij}f_{ij}^2+2\langle \nabla (\Delta f),\nabla
 f\rangle+2R_{ij}f_if_j-\alpha (\Delta f)_t\\
  &-2\alpha
  h_{ij}f_{ij}-2\alpha h_{ik;i}f_k+\alpha H_if_i-\alpha \Delta q\bigg]\\
  &=t(2\sum_{ij}f_{ij}^2-2\alpha h_{ij}f_{ij})+
  2t\langle \nabla (f_t+q-|\nabla f|^2),\nabla
 f\rangle\\
 &- \alpha t(f_{t}+q-|\nabla f|^2)_t -\alpha t(2 h_{ki;k}-H_i)f_i+
 2tR_{ij}f_if_j-\alpha t\Delta q  \\
 &=t(2\sum_{ij}f_{ij}^2-2\alpha h_{ij}f_{ij})+
  2t\langle \nabla (f_t+q-|\nabla f|^2),\nabla
 f\rangle\\
 &+ \alpha t( \frac1\alpha \frac Ft +(1-\frac1 \alpha)|\nabla f|^2)_t
 -\alpha t(2 h_{ki;k}-H_i)f_i+
 2tR_{ij}f_if_j-\alpha t\Delta q  \\
 &=t(2\sum_{ij}f_{ij}^2-2\alpha h_{ij}f_{ij})+
  2t\langle \nabla (\alpha f_t+q-|\nabla f|^2),\nabla
 f\rangle\\
 &+F_t-\frac Ft-2t(\alpha-1)h_{ij}f_if_j-
 \alpha t(2 h_{ki;k}-H_i)f_i+
 2tR_{ij}f_if_j-\alpha t\Delta q\\
 &=t(2\sum_{ij}f_{ij}^2-2\alpha h_{ij}f_{ij})-
  2t\langle \nabla (\frac Ft+(\alpha -1)q),\nabla
 f\rangle\\
 &+F_t-\frac Ft-2t(\alpha-1)h_{ij}f_if_j-
 \alpha t(2 h_{ki;k}-H_i)f_i+
 2tR_{ij}f_if_j-\alpha t\Delta q\\
 &=t(2\sum_{ij}f_{ij}^2-2\alpha h_{ij}f_{ij})-
  2 \langle \nabla F,\nabla f\rangle+F_t-\frac Ft\\
 &-2t(\alpha-1)h_{ij}f_if_j-
 \alpha t(2 h_{ki;k}-H_i)f_i+
 2tR_{ij}f_if_j-\alpha t\Delta q\\
 &-2t(\alpha-1)q_if_i. 
\end{split}
\end{equation}

By \cite{Sh2}, there is a smooth function $\rho$ such that
\begin{equation*}
\left\{\begin{array}{rcl} \frac{1}{C_1}\rho(x)&\leq& r_0(x)\leq
C_1\rho(x)\\|\nabla^0\rho|&\leq& C_1\\|\nabla^0\nabla^0\rho|&\leq&
C_1
\end{array}\right.
\end{equation*}
where $C_1$ is a constant depending on $n$ and the  bound of
$|Rm|$ of $g(0)$. Here $r_0(x)$ is the distance with respect to
$g(0)$ from a fixed point $o$. By the assumption that $|h|$ and
$|\nabla h|$ are uniformly bounded on space-time, we have

\begin{equation*}
\left\{\begin{array}{rcl}\frac{1}{C_2}\rho(x)&\leq& r_t(x)\leq
C_2\rho(x)\\ |\nabla\rho|&\leq& C_2\\|\nabla^2\rho|&\leq& C_2
\end{array}\right.
\end{equation*}
where $C_2$ depends on $C_1$, $T$ and the uniformly upper bound of
$|h|$ and $|\nabla h|$. Here $r_t(x)$ is the distance with respect
to $g(t)$ from $o$.

Let $\eta\in C^{\infty}([0,\infty)$ be such that $0\le \eta\le1$,
$\eta=1$ on $[0,1]$ and $\eta=0$ on $[2,\infty)$, $\eta>0$ on
$[0,2)$, $0\ge \eta'/\eta^\frac12\ge -C_3$ and
 $\eta''\geq -C_3$ on $[0,\infty)$
where $C_3$ is some positive absolute constant. For any $R>0$, let
$\phi=\eta(\rho/R)$.  Suppose at the point $(x_0,t_0)$  where
$\phi F$ attains positive maximum, $0<t_0\le T$. Then at
$(x_0,t_0)$, we have $\phi F_t\ge0$, $F\nabla \phi+\phi\nabla F
=0$ and $\Delta(\phi f)\le0$. Hence at $(x_0,t_0)$:
\begin{equation}\label{max1}
\begin{split}
0&\ge \Delta(\phi F)\\
&=\phi\Delta F+2\langle\nabla \phi,\nabla F\ra+F\Delta \phi\\
&=\phi\Delta F-2F\frac{|\nabla \phi|^2}\phi+F\Delta \phi\\
&\ge \phi\Delta F-C_4F\lf(R^{-1}+R^{-2}\ri)\\
&\ge t_0\phi(2f_{ij}^2-2\alpha h_{ij}f_{ij})-2\phi  \langle \nabla
F,\nabla f\rangle +F_t -\phi\frac F{t_0} \\
&-C_5t_0\phi|\nabla f|^2 -C_5t_0\phi-C_4F\lf(R^{-1}+R^{-2}\ri)\\
&\ge t_0\phi(2f_{ij}^2-2\alpha h_{ij}f_{ij})+2F  \langle \nabla
\phi,\nabla f\rangle  -\phi\frac F{t_0} +\\
&-C_5t\phi|\nabla f|^2 -C_5t\phi-C_4F\lf(R^{-1}+R^{-2}\ri)\\
&\ge t_0\phi(2f_{ij}^2-2\alpha
h_{ij}f_{ij})-C_6F\phi^{\frac12}R^{-1}|\nabla f|  -\phi\frac F{t_0}  \\
&-C_5t_0\phi|\nabla f|^2 -C_5t_0\phi-C_4F\lf(R^{-1}+R^{-2}\ri)\\
&\ge t_0\phi\cdot \frac{2n}{n+\e} f_{ij}^2-C_6F\phi^{\frac12}R^{-1}|\nabla f|  -\phi\frac F{t_0}  \\
&-C_5t_0\phi|\nabla f|^2 -C_7t_0\phi-C_4F\lf(R^{-1}+R^{-2}\ri)\\
&\ge t_0\phi\frac{2}{n+\e}(\Delta f)^2 -C_6F\phi^{\frac12}R^{-1}|\nabla f|  -\phi\frac F{t_0}  \\
&-C_5t_0\phi|\nabla f|^2 -C_7t_0\phi-C_4F\lf(R^{-1}+R^{-2}\ri)\\
&=t_0\phi\frac{2}{n+\e}(|\nabla f|^2-f_t-q)^2-C_6F\phi^{\frac12}R^{-1}|\nabla f|  -\phi\frac F{t_0}  \\
&-C_5t_0\phi|\nabla f|^2 -C_7t_0\phi-C_4F\lf(R^{-1}+R^{-2}\ri).
\end{split}
\end{equation}
Here and below $C_i$'s are constants depending  only on  bounds of
$h$, $|\nabla h|$, $| Rm|$ of $g(0)$, $|\Delta q|$,   $\alpha$,
$\e$, $n$ and $T$.

Multiply both sides by $t_0\phi$, and let $A=\phi|\nabla f|^2$,
$B=\phi(f_t+q)$ so that $\phi F=t_0(A-\alpha B)$. Then for any
$\delta>0$ and $\tau>0$:
\begin{equation}\label{max2}
   \begin{split}
 0&\ge \frac{2t_0^2}{n+\e}(A-B)^2-C_6t_0 \phi F R^{-1}A^\frac12
 -\phi^2 F  -C_5t_0^2\phi A -C_7t_0^2\phi^2\\
 &\qquad-C_4t_0\phi F\lf(R^{-1}+R^{-2}\ri)\\
&\ge \phi
F\lf[-C_4t_0\lf(R^{-1}+R^{-2}\ri)-1\ri]\\
&\qquad+\frac{2t^2_0}{n+\e}\lf[(A-B)^2-
\frac{C_6}2R^{-1}(A-\alpha
B)A^\frac12-C_5A\ri]-C_7t^2_0\\
&\ge \phi
F\lf[-C_4t_0\lf(R^{-1}+R^{-2}\ri)-1\ri]+\frac{2t^2_0}{n+\e}\bigg[(A-B)^2 \\
&\qquad -
\frac{C_6}{2\delta}R^{-2}\delta^{-1}(A-\alpha
B)-\delta(A-\alpha
B)A -C_5A\bigg]-C_7t^2_0\\
&\ge \phi
F\lf[-C_8t_0(1+\delta^{-1})\lf(R^{-1}+R^{-2}\ri)-1\ri]\\
&\qquad+\frac{2t_0^2}{n+\e}\lf[(A-B)^2
- \delta(A-\alpha B)A -\tau A^2\ri]-C_{10}(1+\tau^{-1})t^2_0
\end{split}
\end{equation}
Now for $\sigma>0$:
\begin{equation}\label{max3}
\begin{split}
   (A&-B)^2  - \delta(A-\alpha B)A -\tau A^2 \\
   &=(A-\alpha B)^2+2(\alpha-1)(A-\alpha B)B+(\alpha-1)^2B^2-
   \delta(A-\alpha B)A - \tau A^2 \\
   &=(1-\sigma)(A-\alpha B)^2+ (\sigma-\delta-\tau)A^2+
   [-2\sigma\alpha+2 (\alpha-1)+\delta
   \alpha]AB\\
   &+(\sigma\alpha^2+1-\alpha^2)B^2
\end{split}
\end{equation}
First choose $\sigma$ such that $\sigma\alpha^2+1-\alpha^2=0$,
that is:
$$
\sigma=\frac{\alpha^2-1}{\alpha^2}.
$$
Then $0<\sigma<1$. Next choose $\delta$ so that $-2\sigma\alpha+2
(\alpha-1)+\delta \alpha=0$, that is:
$$
\delta=\frac2\alpha(
\sigma\alpha-\alpha+1)=\frac2\alpha(\frac{\alpha^2-1}{\alpha}-\alpha+1)
=\frac{2(\alpha-1)}{\alpha^2}>0.
$$
Then
$$
\eta-\delta=\frac{(\alpha-1)^2}{\alpha^2}>0.
$$
Then we can choose $\tau=\frac12(\sigma-\delta)$. Then
$$
(A -B)^2  - \delta(A-\alpha B)A -\tau A\ge (1-\sigma)(A-\alpha
B)^2.
$$
Note that $\sigma$, $\delta$ and $\tau$ depend only on $\alpha>1$.
Put this back to (\ref{max2}), we have
\begin{equation}\label{max4}
  \begin{split}
 0&\ge \phi
F\lf[-C_8t_0(1+\delta^{-1})\lf(R^{-1}+R^{-2}\ri)-1\ri]+\frac{2t_0^2}{n+\e}
\lf[(1-\sigma)(A-\alpha B)^2
 \ri]\\
 &\qquad-C_{10}(1+\tau^{-1}) t^2_0\\
 &=\phi F\lf[-C_8t_0(1+\delta^{-1})\lf(R^{-1}+R^{-2}\ri)-1\ri]+
 \frac2{n+\e}\alpha^{-2}( \phi F)^2\\
 &\qquad +\frac{2t_0^2}{n+\e}\lf[(A-B)^2
- \delta(A-\alpha B)A -\tau A^2\ri]-C_{10}(1+\tau^{-1}) t^2_0
\end{split}
\end{equation}
So
$$
\phi F\le C_{11}t_0\lf(R^{-2}+R^{-1}+1\ri)
 +\frac{n+\e}{2}\alpha^2
$$
on $M\times [0,T]$. Since $t_0\le T$,    we see that if $r(x)\le
\frac12 C_2R$, then
$$
F(x,T)\le C_{11}T\lf(R^{-2}+R^{-1}+1\ri)
 +\frac{n+\e}{2}\alpha^2.
$$
Since we can take any $t\in (0,T]$ to be our $T$, the result will
follow by letting $R\to\infty$.
\end{proof}

\begin{cor}\label{gradient-local}
Same assumptions as in the lemma, the following local version of
gradient estimate is true:
\begin{equation*}
\begin{split}
\sup_{B_p(C_2R/2)}&\Big(\frac{|\nabla
u|^2}{u^2}-\alpha\frac{u_t}{u}-\alpha q\Big)\\
&\le \frac{(n+\epsilon)\alpha^2}{2t}+C_{11}\lf(R^{-2}+R^{-1}+1\ri)
\end{split}
\end{equation*}
where the constants $C_i$ are as in the proof of the lemma.
\end{cor}
\begin{rem}\label{gradientlocal1} Note that the constants in the
local gradient estimates depend only on local data and the local
behavior of the function $\rho$
\end{rem}
\begin{cor}\label{cor-harnack}
Let the $u$ be a positive solution of equation
(\ref{eqn-schordinger}). Then, for any $\alpha>1$ and
$\epsilon>0$, there are $C_1>0$ depending on $T$ and the upper
bound of $|h|$, and $C_2>0$ depending on $\alpha,\epsilon, n, T$,
the  upper bounds of $|h|,|\nabla^th|,|q|,|\nabla^tq|,|\Delta^t
q|,|Rc^t|$ and the curvature bound of the initial metric, such
that
\begin{equation*}
u(x_1,t_1)\leq
u(x_2,t_2)\Big(\frac{t_2}{t_1}\Big)^\frac{(n+\epsilon)\alpha}{2}
\exp\left({\frac{C_1\alpha
r^2(x_1,x_2)}{t_2-t_1}+C_2(t_2-t_1)}\right).
\end{equation*}
for any $x_1,x_2\in M$ and $0<t_1<t_2\leq T$.
\end{cor}
\begin{proof}  Let $(x_1,t_1)$ and $(x_2,t_2)$ be two points in
$M\times (0,T]$ with $t_1<t_2$. Let $\gamma(s)$ be a minimal
geodesic joining $x_1$ to $x_2$ with respect to the initial
metric. Let $l=r(x_1,x_2)$. Let $t(s)$ be an affine function such
that $t(0)=t_1$ and $t(l)=t_2$.Then
\begin{eqnarray*}
\log\frac{u(x_2,t_2)}{u(x_1,t_1)}&=&\int_{0}^{l}\frac{d}{ds}\log
u(\gamma(s),t(s))ds\\
&=&\int_{0}^{l}\frac{\la\nabla^t
u,\gamma'\ra}{u}+t'\frac{u_t}{u}ds\\
&\geq&\int_{0}^{l}-C_1\frac{|\nabla^t
u|}{u}+t'\Big(\frac{1}{\alpha}\frac{|\nabla^tu|^2}{u^2}-\frac{C_2}{\alpha}-\frac{(n+\epsilon)\alpha}{2t(s)}-q\Big)ds\\
&\geq&-\frac{\alpha C_1
l^2}{4(t_2-t_1)}-C_3(t_2-t_1)-\frac{n+\epsilon}{2}\log\frac{t_2}{t_1}.
\end{eqnarray*}
where $C_1$ depends only on the  upper bound of $|h|$ and $T$,
$C_3$ depends only on $\alpha,\epsilon, n, T$, the   upper bounds
of $|h|$, $|\nabla^th|$, $|q|$, $|\nabla^tq|$, $|\Delta^t
q|,|Rc^t|$ and the bound of the sectional curvature of the initial
metric. This completes the proof of the Corollary.
\end{proof}

\begin{cor}\label{cor-mean-value-inequlity}
Let $u$ be a positive solution of equation
(\ref{eqn-schordinger}). Then there is a positive constant $C$
depending only on $n, T$,  the   upper bounds of
$|h|,|\nabla^th|,|\nabla^tq|,|\Delta^t q|,|Rc^t|$ and the
curvature bound of $g(0)$, such that for any $x\in M$ and
$0<s<t\leq T$,
\begin{equation*}
u(x,s)\leq
\frac{C}{V_{x}(\sqrt{t-s})}\Big(\frac{t}{s}\Big)^{n+1}\int_{B_{x}(\sqrt{t-s})}u(y,t)dV(y).
\end{equation*}
\end{cor}
\begin{proof} By Corollary \ref{cor-harnack},
\begin{equation*}
u(x,s)\leq C\Big(\frac{t}{s}\Big)^{n+1} u(y,t)
\end{equation*}
for any $y\in B_x(\sqrt{t-s})$, by choosing $\epsilon=1,\alpha=2$,
where $C$ depends on $n, T$,  the   upper bounds of
$|h|,|\nabla^th|,|\nabla^tq|,|\Delta^t q|,|Rc^t|$ and the
curvature bound of $g(0)$. Integrating on the both sides with
respect to $dV(y)$ on $B_x(\sqrt{t-s})$, the result  follows.
\end{proof}
\begin{rem} Since $g(t)$ is uniformly equivalent to $g(0)$, by
volume comparison, we can see that in the corollary, the geodesic
ball and its volume can be chosen with respect to any $g(t)$,
perhaps with a different constant. \end{rem}

\section{Upper and lower estimates of the  fundamental solutions}

In the following, we will apply the last three sections to get
upper and lower estimates for the fundamental solutions of the
equation (\ref{eqn-schordinger}). We always assume
 (\textbf{A1})--(\textbf{A3}) in \S2 are true.

Let $\mathcal{Z}(x,t;y,s)$, $0\le s<t\le T$ be the fundamental
solution of equation (\ref{eqn-schordinger}):
$$
\frac{\p }{\p t} u-\Delta^t u+qu=0.
$$
That is to say:
\begin{equation}\label{fundamental1}
\left\{\begin{array}{l}\frac{\p }{\p t}
\mathcal{Z}(x,t;y,s)-\Delta^t_x
\mathcal{Z}(x,t;y,s)+q(x)\mathcal{Z}(x,t;y,s)=0\\
\lim_{t\to s}\mathcal{Z}(x,t;y,s)=\delta_y.\end{array}\right.
\end{equation}
The fundamental solution exists and is positive, see for example
\cite{Guenther}.

Then $\mathcal{Z}(x,t;y,s)$ is the fundamental solution of the
conjugate equation. That is:
\begin{equation}\label{fundamental2}
\left\{\begin{array}{l}-\frac{\p }{\p s}
\mathcal{Z}(x,t;y,s)-\Delta^s_y
\mathcal{Z}(x,t;y,s)+(q(y)-H(y))\mathcal{Z}(x,t;y,s)=0\\
\lim_{s\to t}\mathcal{Z}(x,t;y,s)=\delta_x.\end{array}\right.
\end{equation}

The fundamental solution $\mathcal{Z}(x,t;y,s)$ can be obtained as
follows.

Let
$\Omega_1\subset\subset\Omega_2\subset\subset\cdots\subset\subset
M$ be an exhaustion of   relatively compact domains with smooth
boundary in $M$. Let $\mathcal{Z}_k(x,t;y,s)$ be the corresponding
fundamental solution on $\Omega_k$ with zero Dirichlet boundary
condition.  Then $\mathcal{Z}_k$ is an increasing sequence by
maximum principle and $\mathcal{Z}$ is the limit of
$\mathcal{Z}_k$ as $k\to\infty$. Moreover, we have
\begin{equation*}
\mathcal{Z}_k(\cdot,\cdot;y,s)\to \mathcal{Z}(\cdot,\cdot;y,s)
\end{equation*}
uniformly on any compact subset of $M\times(s,T]$ up to any
derivatives, and
\begin{equation*}
\mathcal{Z}_k(x,t;\cdot,\cdot)\to \mathcal{Z}(x,t;\cdot,\cdot)
\end{equation*}
uniformly on any compact subset of $M\times[0,t)$ up to any
derivatives.

\begin{lem}\label{lm-int-bound}
There is a positive constant  $C$ depending only on $T$ and the
 upper bounds of $|q|$ and $|H|$, such that
\begin{equation}\label{integralbound1}
\int_{M}\mathcal{Z}(x,t;y,s)dV_t(x)\leq C
\end{equation}
for any $0<s<t\leq T$. Moreover, if $q=H$, then
\begin{equation}\label{integralbound2}
\int_{M}\mathcal{Z}(x,t;y,s)dV_t(x)=1
\end{equation}
for any $0<s<t\leq T$.
\end{lem}
\begin{proof} With the above notations, let
$$
I_k(t)=  \int_{\Omega_k}\mathcal{Z}_k(x,t;y,s)dV_t(x).
$$
Then
\begin{eqnarray*}
& &\frac{d}{dt}I_k(t)\\
&= &\frac{d}{dt}\int_{\Omega_k}\mathcal{Z}_k(x,t;y,s)dV_t(x)\\
&=&\int_{\Omega_k}(\Delta^t_{x}\mathcal{Z}_k-q\mathcal{Z}_k)dV_t(x)+\int_{\Omega_k}H\mathcal{Z}_kdV_t(x)\\
&=&\int_{\Omega_k}(H-q)\mathcal{Z}_kdV_t(x)+\int_{\partial \Omega_k}\frac{\partial Z_k}{\partial {\vec n}_t}dS_t(x)\\
&\leq&\int_{\Omega_k}(H-q)\mathcal{Z}_kdV_t(x)\\
&\leq& C_1I_k(t)
\end{eqnarray*}
since that $\mathcal{Z}_k\geq 0$ on $\Omega_k\times (s,T] $ and it
is $0$ on $\partial \Omega_k\times (s,T]$, where $C_1$ depends on
the uniform upper bounds of $|q|$ and $|H|$. So
\begin{equation*}
\frac{d}{dt}\log I_k(t)\leq C_1.
\end{equation*}
Note that $I_k(s)=1$. Hence
\begin{equation*}
I_k(t)\leq e^{C_1(t-s)}\leq e^{C_1T}.
\end{equation*}
By letting $k\to\infty$, we get the first inequality
(\ref{integralbound1}).

Suppose $q=H$. Let $\phi=\eta(\rho/R)$ be the same as in the proof
of Lemma \ref{lm-harnack}. For any $t_1,t_2$ with $s<t_1<t_2\leq
T$, we have
\begin{eqnarray*}
&&\Big|\int_{M}\phi\mathcal{Z}dV_{t_2}(x)-\int_{M}\phi\mathcal{Z}dV_{t_1}(x)\Big|\\
&=&\Big|\int_{t_1}^{t_2}\int_{M}\phi
(Z_t+H\mathcal{Z})dV_t(x)dt\Big|\\
&=&\Big|\int_{t_1}^{t_2}\int_{M}\phi \Delta^t \mathcal{Z}dV_t(x)dt\Big|\\
&=&\Big|\int_{t_1}^{t_2}\int_{M}\mathcal{Z}\Delta^t \phi
dV_t(x)dt\Big|\\
&\leq& \frac{C_2(t_2-t_1)}{R}\cdot\max_{t_1\leq t\leq
t_2}\int_{M}\mathcal{Z}dV_t\\
&\leq&\frac{C_3}{R}
\end{eqnarray*}
where $C_2,C_3$ are independent of $R$ and (\ref{integralbound1})
has been used. Let $R\to\infty$, we get
\begin{equation*}
\int_{M}\mathcal{Z}dV_{t_2}(x)=\int_{M}\mathcal{Z}dV_{t_1}(x)
\end{equation*}
for any $s<t_1<t_2\leq T$. Note that
\begin{equation*}
\lim_{t\to s^+}\int_{M}\mathcal{Z}dV_{t}(x)=1.
\end{equation*}
The result follows.
\end{proof}
\begin{cor}\label{cor-int-bound}
There is a positive constant $C$ depending on $T$ and the
  upper bounds of $|q|$ and $|H|$, such that
\begin{equation*}
\int_{M}\mathcal{Z}(x,t;y,s)dV_s(y)\leq C
\end{equation*}
for any $s\in [0,t)$.
\end{cor}
\begin{proof} Since $\mathcal{Z}(x,t;y,s)$ is also the fundamental solution of the conjugate
equation, the proof is similar to the proof of Lemma
\ref{lm-int-bound}.
\end{proof}

\begin{lem}\label{prop-p-estimate}
There is a positive constant $C$ depending only on $T,n,$ the
lower bound of the Ricci curvature of the initial metric and the
 upper bounds of $|q|$ and $|h|$, such that
\begin{eqnarray*}
\mathcal{Z}(x,t;y,s)&\leq& \frac{C}{V_x(\sqrt{t-s})}\ \ \mbox{and}\\
\mathcal{Z}(x,t;y,s)&\leq& \frac{C}{V_y(\sqrt{t-s})}.
\end{eqnarray*}
\end{lem}
\begin{proof} Apply the mean value inequality Lemma
\ref{lm-mean-value-inequality} to

 $u(y,s)=\mathcal{Z}(x,t;y,t-s)$
with $r=\frac{\sqrt{s}}{2}$, we get
\begin{eqnarray*}
\mathcal{Z}(x,t;y,t-s)&=&u(y,s)\\
&\leq&\frac{C_1e^{A_1s+B_1r}}{r^2V_{y}(r)}\int_{0}^{s}\int_{M}udV_tds\\
&\leq&\frac{C_2e^{A_1s+B_2\sqrt s}}{V_{y}(\sqrt{s}/2)}\\
&\leq&\frac{C_3e^{A_1s+B_3\sqrt s}}{V_y(\sqrt s)}\\
&\leq&\frac{C_3e^{A_1T+B_3\sqrt T}}{V_y(\sqrt s)}
\end{eqnarray*}
where in the last but second inequality we have  used  volume
comparison. Here $C_1, C_2, C_3$ depend  only on $n$, $T$, and the
  upper bounds of $q$ and $|h|$, $A_1$ depends only the
  upper bounds of $|q|$ and $|H|$, and $B_1$ depends only on
$n$ and the lower bound of the Ricci curvature of $g(0)$. So, we
get the second inequality in the lemma.

Applying similar method to $u(x,t)=\mathcal{Z}(x,t+s;y,s)$ with
$r=\frac{\sqrt{s}}{2}$ will get the first inequality.
\end{proof}

\begin{lem}\label{prop-weighted-L2}
There are some positive constants $C$ and $D$ with $C$ depending
only on $T,n,$ the lower bound of the Ricci curvature of the
initial metric and the  upper bounds of $|q|$ and $|h|$, and $D$
depending only on $T$ and the   upper bound of $|h|$, such that
for $0\le s<t\le T$,
\begin{equation*}
\int_M\mathcal{Z}^2(x,t;y,s)e^{\frac{r^2(x,y)}{D(t-s)}}dV_t(x)\leq
\frac{C}{V_y(\sqrt{t-s})}\ \ \mbox{and}
\end{equation*}
\begin{equation*}
\int_M\mathcal{Z}^2(x,t;y,s)e^{\frac{r^2(x,y)}{D(t-s)}}dV_s(y)\leq
\frac{C}{V_x(\sqrt{t-s})}.
\end{equation*}
\end{lem}
\begin{proof} We just prove the first inequality. The proof of the second
one is similar.

By Lemma \ref{prop-p-estimate} and the fact that $\mathcal{Z}_k$
increasing to $\mathcal{Z}$,
\begin{eqnarray*}
\int_{\Omega_k}\mathcal{Z}^2_{k}(x,t;y,s)dV_t(x)&\leq&
\frac{C_1}{V_y(\sqrt{t-s})}\int_{\Omega_k}\mathcal{Z}_k(x,t;y,s)dV_t(x)\\
 &\leq&\frac{C_2}{V_y(\sqrt{t-s})}
\end{eqnarray*}
for $0\le s<t\le T$, where $C_1$ and $C_2$ depends on $T,n,$ the
lower bound of the Ricci curvature of the initial metric and the
uniformly upper bounds of $|q|$ and $|h|$. Now fix $t>s$ and
consider the function $u(x,\tau)=Z_k(x,\tau+s;y,s)$, $0<
\tau\leq t-s$.

Let
\begin{equation*}
f(\tau)=V_y(\sqrt{\tau}).
\end{equation*}
Then for $0<\tau_1<\tau_2\le T$,
\begin{equation*}
\frac{f(\tau_1)}{f(\tau_1/4)}=\frac{V_y(\sqrt {\tau_1})}{V_y(\sqrt
{\tau_1}/2)}\leq\frac{V_{k}(\sqrt {\tau_1})}{V_k(\sqrt
{\tau_1}/2)}\leq\frac{V_k(\sqrt T)}{V_k(\sqrt T/2)}\le A
\frac{f(\tau_2)}{f(\tau_2/4)},
\end{equation*}
where $V_k(r)$ denotes the volume of the ball of radius $r$ in the
space form with Ricci curvature $-k$ ($-k$ is the lower bound of
the Ricci curvature of the initial metric) and $A=\frac{V_k(\sqrt
T)}{V_k(\sqrt T/2)}$. So, $f$ is regular with the constants $A$
and $\gamma=4$. By Lemma \ref{lem-integral-estimate},
\begin{equation*}
\int_{\Omega_k}\mathcal{Z}^2_{k}(x,t;y,s)e^\frac{r^2(x,y)}{D(t-s)}dV_t(x)\leq
\frac{C}{V_y(\sqrt{t-s})}
\end{equation*}
where $D$ depends on $T$ and the uniformly upper bound of $|h|$,
and $C$ depends on $T,n,$ the lower bound of the Ricci curvature
of the initial metric and the uniformly upper bounds of $|q|$ and
$|h|$.

By taking limit, we get the first inequality.

\end{proof}

\begin{thm}\label{prop-gussian-upper-bound}
There exist positive constants  $C$ and $D$ with $C$ depending
only on $T,n,$ the lower bound of the Ricci curvature of the
initial metric and the   upper bounds of $|q|$ and $|h|$, and $D$
depending only on $T$ and the   upper bound of $|h|$, such that
for $0\le s<t\le T$,
\begin{equation*}
\mathcal{Z}(x,t;y,s)\leq
\frac{C}{V^{\frac{1}{2}}_x(\sqrt{t-s})V^{\frac{1}{2}}_y(\sqrt{t-s})}\times
e^{-\frac{r^2(x,y)}{D(t-s)}}.
\end{equation*}
\end{thm}
\begin{proof} By the triangle inequality, we have
\begin{equation*}
r^2(x,\zeta)+r^{2}(\zeta,y)-\frac{r^2(x,y)}{2}\geq 0.
\end{equation*}
Let $D$ be   as in Lemma \ref{prop-weighted-L2} and let
$\tau=(s+t)/2$.  Then by the semigroup property and Lemma \ref{prop-weighted-L2},
\begin{eqnarray*}
& &\mathcal{Z}(x,t;y,s)\\
&=&\int_{M}\mathcal{Z}(x,t;\zeta,\tau)\mathcal{Z}(\zeta,\tau;y,s)dV_{\tau}(\zeta)\\
&\leq&\int_{M}\mathcal{Z}(x,t;\zeta,\tau)\mathcal{Z}(\zeta,\tau;y,s)
e^{\frac{r^{2}(x,\zeta)}{2D(t-s)}+\frac{r^{2}(\zeta,y)}{2D(t-s)}
-\frac{r^{2}(x,y)}{4D(t-s)}}dV_{\tau}(\zeta)\\
&\leq&e^{-\frac{r^{2}(x,y)}{4D(t-s)}}
\Big(\int_{M}\mathcal{Z}^2(x,t;\zeta,\tau)e^{\frac{r^{2}(x,\zeta)}{D(t-s)}}dV_{\tau}(\zeta)\Big)^{1/2}\\
&&\times\Big(\int_{M}\mathcal{Z}^2(\zeta,\tau;y,s)e^{\frac{r^{2}(\zeta,y)}
{D(t-s) }}dV_{\tau}(\zeta)\Big)^{1/2}\\
&\leq&\frac{C}{V^{\frac{1}{2}}_x(\sqrt{t-s})V^{\frac{1}{2}}_y(\sqrt{t-s})}\times
e^{-\frac{r^2(x,y)}{4D(t-s)}}.
\end{eqnarray*}
\end{proof}

\begin{cor}\label{cor-upper-bound}
There exist positive constants $C$ and $D$ with $C$ depending only
$T,n,$ the lower bound of the Ricci curvature of the initial
metric and the  upper bounds of $|q|$ and $|h|$, and $D$ depending
only on $T$ and the   upper bound of $|h|$, such that
\begin{equation*}
\mathcal{Z}(x,t;y,s)\leq
\frac{C}{V_x(\sqrt{t-s})}e^{-\frac{r^2(x,y)}{D(t-s)}}\ \
\mbox{and}
\end{equation*}
\begin{equation*}
\mathcal{Z}(x,t;y,s)\leq
\frac{C}{V_y(\sqrt{t-s})}e^{-\frac{r^2(x,y)}{D(t-s)}}
\end{equation*}
for any $0<s<t<T$.
\end{cor}
\begin{proof} We just prove the first inequality, the proof of the other
one is similar.

By Proposition \ref{prop-gussian-upper-bound} and volume
comparison,
\begin{eqnarray*}
& &\mathcal{Z}(x,t;y,s)\\
&\leq&
\frac{C_1}{V^{\frac{1}{2}}_x(\sqrt{t-s})V^{\frac{1}{2}}_y(\sqrt{t-s})}\times
e^{-\frac{r^2(x,y)}{D_1(t-s)}}\\
&\leq&\frac{C_1}{V^{\frac{1}{2}}_x(\sqrt{t-s})V^{\frac{1}{2}}_y(\sqrt{t-s})}\times
e^{-\frac{r^2(x,y)}{D_1(t-s)}}\times\frac{V^{\frac{1}{2}}_y(r(x,y)+\sqrt{t-s})}{V^{\frac{1}{2}}_x(\sqrt{t-s})}\\
&\leq&\frac{C_1}{V_x(\sqrt{t-s})}\times
e^{-\frac{r^2(x,y)}{D_1(t-s)}}\times
\frac{V_y(r(x,y)+\sqrt{t-s})}{V_y(\sqrt {t-s})}\\
&\leq&\frac{C_1}{V_x(\sqrt{t-s})}e^{-\frac{r^2(x,y)}{D_1(t-s)}+(n-1)\sqrt
k (r(x,y)+\sqrt{t-s})+n\frac{r(x,y)}{\sqrt{t-s}}}\\
&\leq&\frac{C_1}{V_x(\sqrt{t-s})}e^{-\frac{r^2(x,y)}{D_1(t-s)}+C_2\frac{r(x,y)}{\sqrt{t-s}}}.
\end{eqnarray*}
So, when $\frac{r(x,y)}{\sqrt{t-s}}\geq2C_2D_1$, we have
\begin{equation*}
\mathcal{Z}(x,t;y,s)\leq\frac{C_1}{V_x(\sqrt{t-s})}e^{-\frac{r^2(x,y)}{2D_1(t-s)}}.
\end{equation*}
For those that $\frac{r(x,y)}{\sqrt{t-s}}\leq2C_2D_1$, by
Lemma \ref{prop-p-estimate}, we have
\begin{eqnarray*}
& &\mathcal{Z}(x,t;y,s)\\
&\leq&\frac{C_3}{V_x(\sqrt{t-s})}\\
&=& \frac{C_3}{V_x(\sqrt{t-s})}\times
e^{-\frac{r^2(x,y)}{2D_1(t-s)}}\times
e^{\frac{r^2(x,y)}{2D_1(t-s)}}\\
&\leq&\frac{C_4}{V_x(\sqrt{t-s})}\times
e^{-\frac{r^2(x,y)}{2D_1(t-s)}}.
\end{eqnarray*}
This complete the proof the first inequality.
\end{proof}

Next we want to obtain lower estimates of the fundamental
solution.  We will proceed as in \cite{Davies}.
\begin{lem}\label{lm-int-lower-bound}
There is a positive constant $c$ depending only on $T$ and the
 upper bounds of $|q|$ and $|H|$, such that
\begin{eqnarray*}
\int_{M}\mathcal{Z}(x,t;y,s)dV_t(x)&\geq&c\ \mbox{and}\\
\int_{M}\mathcal{Z}(x,t;y,s)dV_s(y)&\geq& c
\end{eqnarray*}
for any $0<s<t<T$.
\end{lem}
\begin{proof} We just prove the first inequality, and the proof of
the second one is similar. Let $\phi=\eta(\rho/R)$  be the same
function as in the proof of Lemma \ref{lm-harnack}. Then for any
$t_1<t_2$ in $(s,T)$,
\begin{eqnarray*}
&&\frac{d}{dt}\int_{M}\phi \mathcal{Z}dV_t\\
&=&\int_{M}\phi\Delta^t\mathcal{Z}+(H-q)\phi\mathcal{Z}dV_t\\
&\geq&\int_{M}\mathcal{Z}\Delta^t\phi
dV_tdt-C_1\int_{M}\phi\mathcal{Z}dV_t\\
&\geq&-\frac{C_2}{R}-C_1\int_{M}\phi\mathcal{Z}dV_t
\end{eqnarray*}
where $C_1$ depends on the uniformly upper bounds of $|q|$ and
$|H|$, $C_2$ is independent of $R$ and Lemma \ref{lm-int-bound}
has been used in the last inequality. So
\begin{equation*}
\frac{d}{dt}\Big(e^{C_1(t-s)}\int_M\phi \mathcal{Z}dV_t\Big)\geq
-\frac{C_2e^{C_1(t-s)}}{R}\ \mbox{and}
\end{equation*}
\begin{equation*}
\int_M\phi \mathcal{Z}dV_t\geq
e^{-C_1(t-s)}\Big\{e^{C_1(t_1-s)}\int_{M}\phi\mathcal{Z}dV_{t_1}-\frac{C_2(1-e^{C_1(t-s)})}{C_1R}\Big\}
\end{equation*}
for any $T>t>t_1>s$. Let $R\to\infty$ and $t_1\to s^+$, we get
\begin{equation*}
\int_{M}\mathcal{Z}dV_t\geq e^{-C_1(t-s)}\geq e^{-C_1T}
\end{equation*}
for any $t\in(s,T)$.
\end{proof}

\begin{lem}\label{lm-int-lower}
Let $c$ be the constant in Lemma \ref{lm-int-lower-bound}. Then,
there is a constant $A>1$ depending only on $n,T$, the lower bound
of the Ricci curvature of the initial metric and the upper bounds
of $|q|$ and $|h|$, such that
\begin{equation*}
\int_{B_y(A\sqrt{t-s})}\mathcal{Z}(x,t;y,s)dV_t(x)\geq\frac{c}{2}\
\ \mbox{and}
\end{equation*}
\begin{equation*}
\int_{B_x(A\sqrt{t-s})}\mathcal{Z}(x,t;y,s)dV_s(y)\geq\frac{c}{2}
\end{equation*}
for any $0<s<t<T$.
\end{lem}
\begin{proof} we just prove that first inequality, the proof of second one
is similar.

By the second inequality of Corollary \ref{cor-upper-bound},
\begin{eqnarray*}
&&\int_{M\setminus B_y(A\sqrt{t-s})}\mathcal{Z}(x,t;y,s)dV_t(x)\\
&\leq&\frac{C_1}{V_y(\sqrt{t-s})}\int_{M\setminus
B_y(A\sqrt{t-s})}e^{-\frac{r^2(x,y)}{D(t-s)}}dV(x)\\
&=&\frac{C_1}{V_y(\sqrt{t-s})}\int_{A\sqrt{t-s}}^{\infty}e^{-\frac{r^2}{D(t-s)}}dV_y(r)\\
&\leq&C_1\int_{A\sqrt{t-s}}^{\infty}
\frac{V_y(r)}{V_y(\sqrt {t-s})}\times e^{-\frac{r^2}{D(t-s)}}d\Big(\frac{2r}{D(t-s)}\Big)\\
&\leq&\frac{C_1}{D}\int_{A\sqrt{t-s}}^{\infty}\Big(\frac{r}{\sqrt{t-s}}\Big)^ne^{-\frac{r^2}{D(t-s)}+C_2\frac{r}{\sqrt
{t-s}}}d\Big(\frac{r^2}{t-s}\Big)
\end{eqnarray*}
where $C_2$ depends on $n, T$ and the lower bound of the Ricci
curvature of the initial metric.

If we first require that $A\geq 2C_2D$, then
\begin{eqnarray*}
& &\int_{M\setminus B_y(A\sqrt{t-s})}\mathcal{Z}(x,t;y,s)dV_t(x)\\
&\leq&\frac{C_1}{D}\int_{A\sqrt{t-s}}^{\infty}\Big(\frac{r}{\sqrt{t-s}}\Big)^ne^{-\frac{r^2}{2D(t-s)}}d\Big(\frac{r^2}{t-s}\Big)\\
&=&2C_3\int_{2DA^2}^{\infty}x^{\frac{n}{2}}e^{-x^2}dx\\
&\leq&\frac{c}{2}
\end{eqnarray*}
when $A$ is large enough depending only on $n,T$, the lower bound
of the Ricci curvature of the initial metric and the upper bounds
of $|q|$ and $|h|$. This complete the proof of the first
inequality.
\end{proof}
\begin{lem}\label{lm-p-lower-bound}
There is a constant $c>0$ depending on $n, T$ and the uniformly
upper bounds of $|h|,|\nabla^th|,|q|,|\nabla^tq|,|\Delta^t
q|,|Rc^t|$, such that
\begin{equation*}
\mathcal{Z}(x,t;x,s)\geq \frac{c}{V_x(\sqrt {t-s})}\ \mbox{and}
\end{equation*}
\begin{equation*}
\mathcal{Z}(x,t;y,s)\geq \frac{c}{V_y(\sqrt {t-s})}
\end{equation*}
for any $x$, $y$ and $0<s<t<T$
\end{lem}
\begin{proof} We just prove the first inequality, the proof of the second
one is similar.

Let $\tau=\frac{t-s}{2}$. Then, by Corollary \ref{cor-harnack},
\begin{equation*}
c_1e^{-\frac{r^2(x,y)}{c_2\tau}}\mathcal{Z}(y,t-\tau;x,s)\leq
\mathcal{Z}(x,t;x,s)
\end{equation*}
where $c_1>0$ depends only on $n, T$ and the   upper bounds of
$|h|,|\nabla^th|$,$|q|$,
$|\nabla^tq|$,$|\Delta^t q|$,$|Rc^t|$, and $c_2>0$
depends only on $T$ and the   upper bound of $|h|$. This implies
that
\begin{equation*}
c_1\int_{B_x(A\sqrt\tau)}e^{-\frac{r^2(x,y)}{c_2\tau}}\mathcal{Z}(y,t-\tau;x,s)dV(y)\leq
\int_{B_x(A\sqrt t)}\mathcal{Z}(x,t;x,s)dV(y).
\end{equation*}
where $A$ is the same as in Lemma \ref{lm-int-lower}. By Lemma
\ref{lm-int-lower},
\begin{equation*}
\mathcal{Z}(x,t;x,s)\geq
\frac{c_3}{V_x(A\sqrt\tau)}\geq\frac{c_4}{V_x(\sqrt{t-s})}
\end{equation*}
with $c_3$and $c_4$ depending on $n, T$ and the   upper bounds of
$|h|,|\nabla^th|,|q|$,
$|\nabla^tq|,|\Delta^t q|,|Rc^t|$.

This completes the proof of the first inequality.
\end{proof}
\begin{prop}\label{prop-guassian-lower-bound}
There exist positive constants $c$ and $d$ with $c$ depending on
$n, T$ and the   upper bounds of
$|h|,|\nabla^th|,|q|,$ $|\nabla^tq|,|\Delta^t q|,|Rc^t|$, and $d$
depending only on $T$ and the   upper bound of $|h|$, such that
\begin{equation*}
\mathcal{Z}(x,t;y,s)\geq \frac{c}{V_x(\sqrt{t-s})}\times
e^{-\frac{r^2(x,y)}{d(t-s)}}\ \ \mbox{and}
\end{equation*}
\begin{equation*}
\mathcal{Z}(x,t;y,s)\geq \frac{c}{V_y(\sqrt{t-s})}\times
e^{-\frac{r^2(x,y)}{d(t-s)}}.
\end{equation*}
for any $0<s<t<T$.
\end{prop}
\begin{proof} We just prove the second inequality, the proof of the first
one is similar.

Let $\tau=\frac{t-s}{2}$. By Corollary \ref{cor-harnack} and Lemma
\ref{lm-p-lower-bound},
\begin{eqnarray*}
& &\mathcal{Z}(x,t;y,s)\\
&\geq& c_1\mathcal{Z}(y,t-\tau;y,s)e^{-\frac{r^2(x,y)}{c_2\tau}}\\
&\geq& \frac{c_3}{V_y(\sqrt \tau)}e^{-\frac{r^2(x,y)}{c_2\tau}}\\
&\geq& \frac{c_4}{V_y(\sqrt{t-s})}e^{-\frac{r^2(x,y)}{c_2(t-s)}}.\
\end{eqnarray*}
\end{proof}
\begin{cor}\label{cor-lower-bound}
For the same positive constants $c$ and $d$ as in Proposition
\ref{prop-guassian-lower-bound},
\begin{equation*}
\mathcal{Z}(x,t;y,s)\geq
\frac{c}{V^{\frac{1}{2}}_x(\sqrt{t-s})V^\frac{1}{2}_y(\sqrt{t-s})}\times
e^{-\frac{r^2(x,y)}{d(t-s)}}
\end{equation*}
for any $0<s<t<T$.
\end{cor}
\begin{proof} Multiplying the two inequalities in Proposition
\ref{prop-guassian-lower-bound}   the result follows.
\end{proof}

\section{More gradient estimates}

In this section we want to obtain more gradient estimates, which
are generalizations of the gradient estimates in \cite{Zhang} and
\cite{N1} to complete noncompact manifolds. The estimates will be
used in later sections.

Let $u>0$ be a solution of the equation:
\begin{equation}\label{forward}
\Delta^t u-u_t=0
\end{equation}
on $M\times[0,T]$  corresponding to the   Ricci flow
\begin{equation}\label{ricciflow}
   \frac{\p }{\p t}g= -2Ric
\end{equation}
or a solution of
\begin{equation}\label{backward}
   \Delta^{\tau} u-u_\tau-Ru=0
\end{equation}
on $M\times[0,T]$  corresponding to the backward Ricci flow
\begin{equation}\label{backricciflow}
   \frac{\p }{\p \tau}g= 2Ric
\end{equation}
where $R$ is the scalar curvature, $\tau=T-t$.

 Let use make the following assumption:
\begin{enumerate}
    \item [{\bf (a1)}] $g(t)$ is complete and $|\nabla^k Rm|$ are uniformly bounded on
    spacetime by $c_k$ for all $k$.
\end{enumerate}
\begin{lem}\label{gradientbound1} Let $u>0$ be a solution of
(\ref{backward}) or (\ref{forward}) such that $ u \le A$ for all
$t$.
 There is a constant $C$ depending only on those constants $c_k$ for $k=0, 1$ in the
assumption (\textbf{a1}), $n$, $T$ and $A$, such that
$$
 |\nabla u|(x,t)\le C/t
 $$
 for all $x\in M\times(0,T]$.
 \end{lem}
Before we prove the lemma, let us modify a maximum principle in
\cite{Ecker-Huisken}.
\begin{lem}\label{maximumprinciple}
Suppose $g(t)$ is a smooth family of complete metrics defined on
$M$, $0\le t\le T$ with Ricci curvature bounded from below and
$|\frac{\p}{\p t}g|\le C$ on $M\times [0,T]$. Suppose $f(x,t)$ is
a smooth function defined on $M\times [0,T]$ such that
\begin{equation}\label{condition2}
\Delta^t f-\frac{\p}{\p t}f\ge0
\end{equation}
and
\begin{equation}\label{condition1}
\int_0^T\int_M \exp(-ar^2_t(o,x))f^2(x,t)dV_t<\infty
\end{equation}
for some $a>0$. If $f(x,0)\le 0$ for all $x\in M$, then $f\le0$ on
$M\times[0,T]$.
\end{lem}
\begin{proof}
In \cite{Ecker-Huisken}, the condition (\ref{condition1}) with $f$
being replaced by $|\nabla f|$ is assumed. From the their proof,
it is easy to see that the result is still true if $f$ is replaced
by $|\nabla f_+|$, where $f_+=\max\{f,0\}$.

Observe that by (\ref{condition1}), there exists $T_i\uparrow T$
such that $$\int_M
\exp(-ar^2_{T_i}(o,x))f^2(x,T_i)dV_{T_i}<\infty.$$
 Then one can obtain their
condition by using (\ref{condition1}) and (\ref{condition2}) by a
cutoff argument on $[0,T_i]$, perhaps with another $a_i>0$ .
\end{proof}

\begin{proof}[Proof of Lemma \ref{gradientbound1}] In the following
$C_i$'s will  denote constants depending only on the quantities
mentioned in the lemma.

 Let us prove the case that $u$ is a
solution of (\ref{backward}). The other case is similar. Note that
we may consider the interval $[\epsilon,T]$ first, and then let
$\epsilon\to0$. Hence we may assume that $u$ is smooth up to
$\tau=0$. For simplicity, we will write $\tau$ as $t$, and
$\Delta$ instead of $\Delta^t$.

Let us compute the followings in normal coordinates at a point
with respect to $g(t)$
\begin{equation}\label{estimate1}
\begin{split}
(\Delta-\partial_t)|\nabla
u|^2&=2\sum_{ij}u_{ij}^2+2\sum_{i}(\Delta
u)_iu_i+4\sum_{ij}R_{ij}u_iu_j-2\sum_i u_{ti}u_i\\
&=2\sum_{ij}u_{ij}^2+4\sum_{ij}R_{ij}u_iu_j+2\sum_{i}(Ru)_i u_i\\
&\ge 2\sum_{ij}u_{ij}^2 -C_1(|\nabla u|^2+1).
\end{split}
\end{equation}
Here and below, $C_i$ denotes constant depending only on $n,T$ and
$c_0, c_1$. Let $f =\lf(|\nabla u|^2+1\ri)^\frac12$. Then
\begin{equation}\label{estimate2}
\begin{split}
4f^2|\nabla f|^2&=|\nabla f^2|^2\\
&=4\sum_{j}(\sum_{i}u_{ij}u_i)^2\\
&\le 4f^2\sum_{ij}u_{ij}^2.
\end{split}
\end{equation}
\begin{equation}
\begin{split}
 2f(\Delta-\partial_t)f&=  (\Delta-\partial_t)f^2-2|\nabla f|^2\\
 &\ge -C_1f^2
\end{split}
\end{equation}
and so
\begin{equation}\label{estimate2}
(\Delta-\partial_t)(tf)\ge - C_1 f.
\end{equation}
\begin{equation}\label{estimate3}
(\Delta-\partial_t)(\frac12{C_1}  u^2-C_2 t-\frac12C_1 A^2 )=C_1
|\nabla u|^2 +C_1Ru^2+C_2 \ge C_1f^2\ge C_1f
\end{equation}
where $C_2$ is chosen so that $C_1Ru^2+C_2-C_1\ge0$.

 Hence
\begin{equation}\label{estimate4}
(\Delta-\partial_t)Q  \ge 0
\end{equation}
where $$ Q= tf+\frac2{C_1} u^2-C_2 t-\frac12C_1 A^2. $$
 Since
$Q\le 0$, at $t=0$, the result will follow from Lemma
\ref{maximumprinciple}, provided that we can prove:
$$
\int_0^T\int_M\exp(-a^2 r_t^2) Q^2dV_t dt<\infty
$$
for some $a>0$. Here $r_t(x)$ is the distance from a fixed point.
Since the curvature is bounded, by volume comparison, it is
sufficient to prove that
\begin{equation}\label{condition3}
\int_0^T\int_M\exp(-a^2 r_t^2) |\nabla u|^2dV_t dt<\infty
\end{equation}
for some $a>0$. Here $r_t(x)$ is the distance from a fixed point
in $g(t)$. Since we have
$$
\lf(\Delta-\frac{\p}{\p t}\ri)(\exp(-C_3t)u)\ge 0.
$$
Using a cutoff argument and the fact that $u$ is bounded, it is
easy to see that (\ref{condition3}) is true.
\end{proof}
\begin{lem}\label{gradientbound2} With the same notations and
assumptions as before, the following estimates are true:

  Suppose $u>0$ is a solution of (\ref{backward}) or
    (\ref{forward}), and suppose
      $u\le A$. Then $$
t\frac{|\nabla u|^2}u\le C\lf[u\log \frac Au+u\ri]
$$
for some constant $C$ depending only on $T$,  $n$ and  $c_k$,
$0\le k\le 2$, in the assumption (\textbf{a1}).
\end{lem}
\begin{proof} We prove the case that $u$ is a solution of (\ref{backward}) and the
other case can be proved similarly. As in \cite{N2}, for some
suitable positive constants $C_1, C_2, C_3 $  such that if we let
$$
\Phi=\varphi \frac{|\nabla
u|^2}u-\exp(C_1\tau)u\log\lf(\frac{A}u\ri)-C_2\tau u
$$
where $\varphi=\tau/(1+C_3\tau)$, then
$$
\lf(\Delta^t-\frac{\p}{\p \tau}\ri)\Phi\ge0.
$$
 Here and below, $C_i$'s will
denote constants depending on quantities mentioned in the lemma.
Again, we may consider $u$ being a solution in the interval
$[\e,T]$ for $\e>0$ first, and $\tau=0$ corresponding to the
original time $\e$. Let $p\in M$ be a fixed point, then by
Corollary \ref{cor-harnack}
$$
u(p,\frac\e2)\le (C_4+c) u(x,\tau)\exp((C_4+c) r^2(p,x))
$$
for some $c$ depends only on quantities mentioned in the lemma and
$\e$, and $\e\le \tau\le T$. By Lemma \ref{gradientbound1}, and
the fact that $u$ is bounded, we conclude that there is $a>0$ such
that $\exp (-a r_t(p,x)^2)\Phi^2$ is integrable on $M\times[\e,T]$
with respect to $dV_\tau d\tau$. Since $\Phi\le 0$ initially, we
can apply Lemma \ref{maximumprinciple} to conclude that the lemma
is true for bounded and positive solutions of  (\ref{backward}).
\end{proof}

\section{A Li-Yau-Hamilton type differential inequality}
Let $(M^n,g(t))$ be a solution of the Ricci flow:
\begin{equation}\label{ricciflow}
   \frac{\p g}{\p t}=- 2Ric
\end{equation}
on $M\times[0,T]$ for some $T>0$. We always assume that $M^n$ is
noncompact, $g(t)$ is complete and {\bf (a1)} in the previous
section is true. Let $\mathcal{Z}(x,t;y,s)$ with $0\le s<t\le T$
be the fundamental solution of
\begin{equation}\label{heatequation}
  \frac{\p u}{\p t}-\Delta^t u =0.
\end{equation}
Let $p\in M$ be fixed and let $u(x,t)=\mathcal{Z}(p,T;x,t)>0$.
Then $u$ is a solution of the conjugate heat equation

\begin{equation}\label{conjugateheat}
   -\frac{\p u}{\p t} - \Delta^t  u+\mathcal{R}u=0
\end{equation}
on $M\times[0,T]$,  where $\mathcal{R}$ is the scalar curvature
and $ \Delta^t$ is the Laplacian with respect to $  g(t)$. When
there is no confusion, we simply denote $\Delta^t$ as $\Delta$.
  Let $v$ be defined by
\begin{equation}\label{Harnackex}
    v=\lf[\tau(2\Delta f-|\nabla f|^2+R)+(f-n)\ri]u
\end{equation}
where $f$ is defined by $u=e^{-f}/(4\pi\tau)^\frac n2 $ and
$\tau=T-t$ (This notation is adopted throughout this section).

Let $h_0\ge0$ be a smooth function with compact support  and let
$0<t_0<T$. Let $h(x,t)$ be the solution of (\ref{heatequation}) on
$M\times[t_0,T]$ with initial data $h(x,t_0)=h_0(x)$. We want to
prove the following:
\begin{thm}\label{LYHineq} With the above notations and assumption
(a1), we have
\begin{enumerate}
    \item[(i)] For any $t_0<t<T$.
    $hv(\cdot,t)\in L^1(M,  g(t))$,
    \item [(ii)] For any $t_0<t_1<t_2<T$
    $$
    \int_M hvdV_{t_1}\le
    \int_M hvdV_{t_2}.
$$
    \item [(iii)]
    $$
    \limsup_{t\to T^-}\int_M
    hvdV_{t}\le0.
    $$
\end{enumerate}
\end{thm}

\begin{lem}\label{h1} Theorem \ref{LYHineq}(i) is true.
\end{lem}
\begin{proof} For any $T>t>t_0$, by Corollary \ref{cor-upper-bound} and Lemma \ref{gradientbound2},
 there are $a,  C_1>0$
such that for all $x\in M$
\begin{equation}\label{Fundamental1}
   \Big(\frac{|\nabla u|^2}{u}+u+ h+|\nabla h|\Big)(x,t)
   \le C_1\exp(-ar^2(x)).
\end{equation}
Since the curvature is bounded and since $-f=\log u+\frac
n2\log(4\pi\tau)$,   $|\nabla f|^2uh$, $fuh$, $nuh$, $Ruh$ are all
in $L^1$, where $\tau=T-t$. Moreover, since
$$
\Delta f=-\frac{\Delta u}{u}+|\nabla f|^2
$$
in order to prove the lemma, it is sufficient to prove that
$h\Delta u$ is in $L^1$. By Lemma \ref{lm-harnack}, we have for
$T>t>0$,

$$
\frac{|\nabla u|^2}{u^2}-\alpha\frac{u_\tau}u-C_2\le 0
$$
for some $C_2>0$. Hence
\begin{equation}
\begin{split}
\int_M|h\Delta u|&=\int_M |u_\tau-Ru|h\\
&\le \alpha^{-1}\int_M (\alpha u_\tau-\frac{|\nabla u|^2}u +C_2 )h
+C_3\\
&=\int_Mh\Delta u+C_4\\
&<C_5
\end{split}
\end{equation}
for some constants $C_3-C_5$, where have has used
(\ref{Fundamental1}), integration by parts together with a cutoff
argument. This completes the proof of the lemma.
\end{proof}
\begin{rem}\label{h11} From the proof, it is easy to see that for
$0<\tau_1<\tau_2<T$, there is a constant $C$ so that
$$
\int_M(|vh|+|v|)dV_t\le C
$$
for all $\tau_1\le \tau\le\tau_2$. Moreover, (\ref{Fundamental1})
is true for a constant $C$ for all $\tau_1\le \tau\le\tau_2$.
\end{rem}
\begin{lem}\label{h2} Theorem \ref{LYHineq}(ii) is true.
\end{lem}
\begin{proof} By \cite{P1}, we have
\begin{equation}\label{Harnackex1}
    \lf(\frac{\p}{\p \tau}- \Delta
    +R\ri)v=-2\tau\lf|R_{ij}+f_{ij}-\frac1\tau
    g_{ij}\ri|^2u,
\end{equation}
where $\tau=T-t$.  Let $\rho$ and $\phi$ be the functions defined
in the proof of Lemma \ref{lm-harnack}.  Fix $0<t_1<t_2<T$, by
Remark \ref{h11} and (\ref{Fundamental1}) for any $ t_1\le t\le
t_2$,

\begin{equation}\label{monotone1}
\begin{split}
    \frac{d}{d\tau}\int_M \phi hvdV_t &=
    \int_M \phi(v_\tau h +vh_\tau+
    vhR)dV_t \\
        &\le \int_M\phi(h  \Delta v-v
        \Delta h)dV_t\\
        &= \int_M \lf(v \la\nabla \phi,\nabla h\ra-
        h\la\nabla \phi,\nabla v\ra\ri)dV_t\\
        &=\int_M \lf(2v \la\nabla \phi,\nabla h\ra+hv\Delta
        \phi\ri)
        dV_t\\
        &\le \frac{C}{R}
\end{split}
\end{equation}
for some $C>0$ for all $\tau_1\le \tau\le\tau_2$.  By integrating
from $\tau_1$ to $\tau_2$, and letting $R\to\infty$, the result
follows.
\end{proof}

Next, we want to prove Theorem \ref{LYHineq}(iii). We need several
lemmas.

\begin{lem}\label{est1} For any $\alpha>1$ and $\delta, \e>0$, there is a
constant $C(\alpha,\delta,\e)$ which is independent of $t$ such
that, if $\frac{T}{2}<t< T$, then
\begin{equation}\label{taueq02}
 (1-2\alpha\delta)\int_M  \frac{|\nabla u|^2}{u}hdV_t\le
C +\frac{(n+\e)\alpha^2}{2\tau}\int_MuhdV_t
\end{equation}
where $C$ is a constant independent of $t$ and $\tau=T-t$.
\end{lem}
\begin{proof}
Let $\alpha$ and $\e$ be given, then by Lemma \ref{lm-harnack},
$$
\frac{|\nabla u|^2}{u}\le \alpha(u_\tau+Ru)+C_1
+\frac{(n+\e)\alpha^2}{2\tau}= \alpha \Delta
u+C_1+\frac{(n+\e)\alpha^2}{2\tau}.
$$
Here and below, $C_i$'s are positive constant which is independent
of $\tau$.  Let $\rho$ and $\phi$ be as in the proof of Lemma
\ref{lm-harnack}. Then
\begin{equation}\label{taueq03}
   \begin{split}
     \int_M\phi^2 \frac{|\nabla u|^2}{u}hdV_t& \le
     \alpha\int_M\phi^2h \Delta udV_t+C_1\int_M\phi^2 uhdV_t\\
     &+
     \frac{(n+\e)\alpha^2}{2\tau}\int_M\phi^2uhdV_t \\
       & \le \alpha\lf(
       -2\int_M\phi h\la \nabla u,\nabla \phi\ra dV_t
       -\int_M\phi^2\la\nabla u,\nabla h\ra dV_\tau\ri)\\
       &\quad+C_2 +\frac{(n+\e)\alpha^2}{2\tau}\int_MuhdV_t
       \\
       &\le 2\alpha \delta
       \int_M\phi^2 \frac{|\nabla
       u|^2}{u}hdV_t+\frac{C_3\alpha}{\delta R^2}\int_M
       hudV_t\\
       &\quad+\frac{\alpha}{\delta}\int_M\phi^2\frac{|\nabla
       h|^2}{h}u dV_t+C_2 +\frac{(n+\e)\alpha^2}{2\tau}\int_MuhdV_t
   \end{split}
\end{equation}
where $C_i$'s do not depend on $\tau$.  Note that $\int_M h udV_t$
 is bounded by a constant independent of $\tau \in (0, T/2)$
 because $h$ is bounded and
 $$
 \int_M u dV_t
 $$
 is bounded  independent of $\tau$ by Lemma \ref{lm-int-bound}.
Also, $|\nabla h|^2/h$  is bounded independent of $\tau \in (0,
T/2)$ by Lemma \ref{gradientbound2}.  Let $R\to\infty$, we have
\begin{equation}\label{taueq04}
(1-2\alpha\delta)\int_M  \frac{|\nabla u|^2}{u}hdV_t\le
C_1(\alpha,\e,\delta)+\frac{(n+\e)\alpha^2}{2\tau}\int_MuhdV_t.
\end{equation}
\end{proof}
\begin{lem}\label{est2}
For any $\delta>0$,
$$
\int_M(-\Delta u)hdV_t\le 2\delta \int_M  \frac{|\nabla
u|^2}uhdV_t
     +C
     $$
     where $C$ is a constant independent of $t$, provided
     $0<\tau<\frac T2$.
\end{lem}
\begin{proof} Let $\phi$ be then same as before. As in the proof of Lemma \ref{est1},
\begin{equation}\label{taueq05}
     \int_M\phi^2(- \Delta u)hdV_\tau \le
     2\delta \int_M \phi^2\frac{|\nabla u|^2}uhdV_t
     +\frac{C}{\delta R^2}\int_M hu dV_t+C
\end{equation}
where $C$ is independent of $\tau$. Since $h \Delta u$ is in
$L^1(M,  g(t))$ by the proof of Lemma \ref{h1}, let $R\to\infty$,
the result follows.
\end{proof}
\begin{lem}\label{est3}
$$
\limsup_{t\to T^-}\int_M  \tau h(2 \Delta f-|\nabla
f|^2+\mathcal{R})u dV_t\le \frac n2h(x,T).
$$
\end{lem}
\begin{proof} Let $\alpha>1$, $\delta$, $\e>0$ be constants to be chosen later.
By Lemmas \ref{est1} and \ref{est2}:
\begin{equation}\label{taueq06}
     \begin{split}
    \int_M  \tau h(2 &\Delta f-|\nabla
f|^2+R)udV_t    =\int_M \tau h(-2 \Delta u+\frac{|\nabla u|^2}{u}
+
Ru)dV_t\\
         &\le \tau(1+4\delta)\lf[\int_M \frac{|\nabla
         u|^2}{u}hdV_\tau+C_1\ri]+\tau\int_M hu\mathcal{R}dV_t\\
         &\le \tau\frac{1+4\delta} {(1-2\alpha\delta)}
\lf[C_2 +\frac{(n+\e)\alpha^2}{2\tau}\int_MuhdV_\tau\ri]+\tau C_3
     \end{split}
\end{equation}
where $C_1-C_3$ are constants independent of $\tau$, provided
$0<\tau<\frac T2$. Here we have used the fact that $h$ is bounded
and $\int_M udV_t $ is uniformly bounded independent of $\tau$ by
Lemma \ref{lm-int-bound}. Choose $\delta$ small such that
$1-2\alpha\delta>0$.  Then
\begin{equation}\label{taueq07}
\begin{split}
 \limsup_{t\to T^-} \int_M  \tau h(2&\Delta f-|\nabla
f|^2+\mathcal{R})udV_\tau\\
&\le \frac{(n+\e)(1+4\delta)\alpha^2}{(1-2\alpha\delta)}\limsup_{t\to T^-}\int_M uhdV_t \\
&=\frac{(n+\e)(1+4\delta)\alpha^2} {2(1-2\alpha\delta)} h(p,T)
 \end{split}
\end{equation}
where we have used the fact that $u\to \delta_p$ as $t\to T^-$ and
that $h$ is smooth and bounded. Let $\alpha\to1$, $\e,\
\delta\to0$, the result follows.
\end{proof}

\begin{lem}\label{est4}
$$
\limsup_{t\to T^-}\int_M fuhdV_t\le \frac n2 h(x,T)
$$
\end{lem}

\begin{proof} Let $\delta>0$ be fixed and choose $C_1>0$ such that
 $C_1^{-1}\tau^{\frac n2}\le V^t_p(\sqrt
\tau)\leq C_2 \tau^{\frac{n}{2}}$ for any $\tau\in [0,\delta]$
where $\tau=T-t$.

i) We claim that, for any $\epsilon>0$, there is a constant $A>0$,
such that
\begin{equation*}
\D\int_{M\setminus B^{t}_p(A\sqrt{\tau})}fhudV_t\leq \epsilon
\end{equation*}
for any $\tau\in (0,\delta]$.

By   Corollary \ref{cor-upper-bound}, we have
\begin{equation*}
u(x,\tau)\leq
\frac{C_1}{V^{t}_p(\sqrt\tau)}e^{-\frac{r^2_{t}(x,p)}{D\tau}}
\end{equation*}
for some positive constant $C_2$ and $D$. Then

\begin{eqnarray*}
\int_{M\setminus B^{t}_p(A\sqrt\tau)}fhudV_t&=&\int_{M\setminus
B^{t}_p(A\sqrt\tau)}2hu\log\frac{(\frac{1}{4\pi
\tau})^{\frac{n}{4}}}{\sqrt
u}dV_t\\
&\le&\int_{M\setminus
B^{t}_p(A\sqrt\tau)}2hu\cdot\frac{(\frac{1}{4\pi
\tau})^{\frac{n}{4}}}{\sqrt u}dV_t\\
&\leq&\frac{C_3}{\tau^{\frac n4}}\int_{M\setminus
B^{t}_p(A\sqrt\tau)}\sqrt u dV\\
&\leq&\frac{C_4}{\tau^\frac n2}\int_{M\setminus B^{t}_p(
A\sqrt\tau)}e^{-\frac{r_{t}^2(x,p)}{2D\tau}}dV_t(x)\\
&=&\frac{C_4}{\tau^\frac n2}\int_{A\sqrt\tau}^{\infty}e^{-\frac{r^2}{2D\tau}}A_p^t(r)dr \\
&=&\frac{C_4}{\tau^\frac
n2}\lf(\lf[e^{-\frac{r^2}{2D\tau}}V_p^t(r)\ri]_{A\sqrt
\tau}^\infty+\int_{A\sqrt\tau}^{\infty}V_p^{t}(r)e^{-\frac{r^2}{2D\tau}}d(\frac{r^2}{2D\tau})\ri)\\
 &\leq&  C_5 \int_{A\sqrt\tau}^{\infty} \tau^{-\frac n2} \sinh^{n-1} (C_6r)
 e^{-\frac{r^2}{2D\tau}}d\lf(\frac{r^2}{\tau}\ri)\\
&\leq& C_7
\int_{A\sqrt\tau}^{\infty}\lf(\frac{r}{\sqrt\tau}\ri)^{n}e^{C_6r}e^{-\frac{r^2}{2D\tau}}d\lf(\frac{r^2}{\tau}\ri)\\
&=&C_7\int_{A^2}^{\infty}\rho^{\frac{n}{2}}e^{-\frac{\rho}{2D}}
e^{C_6\sqrt{\rho\tau}}d
\rho\\
&\leq&
C_7\int_{A^2}^{\infty}\rho^{\frac{n}{2}}e^{-\frac{\rho}{4D}}d\rho\
\ \
(\mbox{if}\ A\geq 4D C_6\sqrt\delta)\\
&\leq& \epsilon \ \ \ \mbox{(if we choose $A$ large enough),}
\end{eqnarray*}
where $C_1$--$C_7$ are constants independent of $A$ and $\e$. Here
we have used the following facts: $\log x\le x$; $h$ and $u$ are
positive and $h$ is bounded; $V_p^t(r)\le C(n)\sinh^n(C_6r)$ by
volume comparison; $\sinh(C_6r)/r\le Ce^{C_6r}$ for some constant
$C$ depending only on $C_6$.

ii) By the asymptotic behavior of the heat kernel, see
\cite{Guenther} for example, there is an open neighborhood $U$ of
$p$, some positive constants $\tau_0$ and $C_6$ and a positive
function $u_0\in C^{\infty}(U\times [0,\tau_0])$ with
$u_0(p,0)=1$, such that
\begin{eqnarray*}
\Big|u-\frac{1}{(4\pi
\tau)^{\frac{n}{2}}}e^{-\frac{r_t^{2}(x,p)}{4\tau}}u_0(x,\tau)\Big|\leq
C_6\tau^{1-\frac{n}{2}}
\end{eqnarray*}
for any $x\in U$ and $\tau\in(0,\tau_0]$. Hence, for any $x\in
B^{t}_p(A\sqrt\tau)$ when $\tau$ is small,
\begin{eqnarray*}
u&\geq&\frac{1}{(4\pi
\tau)^{\frac{n}{2}}}e^{-\frac{r_t^{2}(x,p)}{4\tau}}u_0(x,\tau)-C_6\tau^{1-\frac{n}{2}}\\
&\geq&\frac{1}{(4\pi
\tau)^{\frac{n}{2}}}e^{-\frac{r_t^{2}(x,p)}{4\tau}}\lf(1-\frac{C_7\tau
}{u_0(x,\tau)}e^{\frac{r^2_{t}(x,p)}{4\tau}}\ri)u_0(x,\tau)\\
&\geq&\frac{1}{(4\pi
\tau)^{\frac{n}{2}}}e^{-\frac{r_t^{2}(x,p)}{4\tau}}(1-C_8\tau)u_0(x,\tau)
\end{eqnarray*}
where all the constants are independent of $\tau$. So
\begin{equation*}
f(x,t)\leq\frac{r_t^{2}(x,p)}{4\tau}-\log(1-C_8\tau)-\log
u_0(x,\tau)
\end{equation*}
for any $x\in B^{t}_p(A\sqrt\tau)$ when $\tau$ is small enough.

Hence
\begin{eqnarray*}
&&\int_{B^{t}_p(A\sqrt\tau)}fhudV_t\\
&\leq&\int_{B^{t}_p(A\sqrt\tau)}\Big(\frac{r_t^{2}(x,p)}{4\tau}-\log(1-C_8\tau)-\log u_0(x,\tau)\Big)hudV_t\\
&\leq&\int_{B^{t}_p(A\sqrt\tau)}\frac{r_t^{2}(x,p)}{4\tau}hudV_t+C_9\tau\\
&=&\int_{B^{t}_p(A\sqrt\tau)}\frac{r_t^{2}(x,p)}{4\tau}\Big(\frac{1}{(4\pi
\tau)^{\frac{n}{2}}}e^{-\frac{r_t^{2}(x,p)}{4\tau}}u_0(x,\tau)+C_6\tau^{1-\frac{n}{2}}\Big)hdV_t+C_9\tau\\
&\leq&\int_{B^{t}_p(A\sqrt\tau)}\frac{r_t^{2}(x,p)}{4\tau}\frac{1}{(4\pi
\tau)^{\frac{n}{2}}}e^{-\frac{r_t^{2}(x,p)}{4\tau}}u_0(x,\tau)hdV_t\\
&+&C_6\int_{B_p^{t}(A\sqrt\tau)}\frac{r^2_{t}(x,p)}{\tau^{\frac{n}{2}}}hdV_t+C_9\tau\\
&=&\frac{1}{(4\pi
\tau)^{\frac{n}{2}}}\int_{0}^{A\sqrt\tau}\frac{r^2}{4\tau}e^{-\frac{r^2}{4\tau}}A^{t}_p(r)\tilde
h(r,t)dr+C_{10}\tau
\end{eqnarray*}
where
\begin{equation*}
\tilde h(r,t)=\frac{1}{A^t_p(r)}\int_{\partial
B^t_{p}(r)}hu_0dS_t.
\end{equation*}
Hence
\begin{eqnarray*}
& &\int_{B^{t}_p(A\sqrt\tau)}fhudV_t\\
&\leq&\frac{1}{2\pi^{\frac{n}{2}}}\int_{0}^{\frac{A^2}{4}}\rho^\frac{n}{2}e^{-\rho}\frac{A^{t}_p(2\sqrt{\tau\rho})}{(2\sqrt{\tau\rho})^{n-1}}\tilde
h(2\sqrt{\tau\rho},t)d\rho+C_{10}\tau
\end{eqnarray*}
where $C_9$ and $C_{10}$ are both independent of $\tau$.

Note that
\begin{equation*}
\frac{A^{t}_p(2\sqrt{\tau\rho})}{(2\sqrt{\tau\rho})^{n-1}}\to\alpha_{n-1}\
\ \mbox{uniformly for}\ \rho\in\lf[0,\frac{A^2}{4}\ri]\ \mbox{as}\
\tau\to0^+,
\end{equation*}
where $\alpha_{n-1}$ means the volume of the the standard sphere
of dimension $n-1$. Moreover
\begin{eqnarray*}
 &
&\tilde
h(2\sqrt{\tau\rho},t)\\&=&\frac{1}{A^t_p(2\sqrt{\tau\rho})}\int_{\partial
B^t_{p}(2\sqrt{\tau\rho})}hu_0dS_t\to h(p,T)u_0(p,0)\\&=&h(p,T)\
\mbox{uniformly for}\ \rho\in\lf[0,\frac{A^2}{4}\ri]\ \mbox{as}\
\tau\to0^+.
\end{eqnarray*}
So,
\begin{eqnarray*}
& &\lim\sup_{\tau\to0^+}\int_{B^{t}_p(A\sqrt\tau)}fhudV_t\\
&\leq&\frac{\alpha_{n-1}h(p,T)}{2\pi^{\frac{n}{2}}}\Gamma(\frac{n}{2}+1)\\
&=&\frac{n\omega_n}{2}\cdot\frac{\Gamma(\frac{n}{2}+1)}{\pi^{\frac{n}{2}}}h(p,T)\\
&=&\frac{n}{2}h(p,T)
\end{eqnarray*}
where $\omega_n$ means the volume of the unit ball in
$\mathbb{R}^n$.

The lemma follows from i) and ii).

\end{proof}

\begin{proof}[Proof of Theorem \ref{LYHineq}(iii)] The result follows from Lemmas \ref{est3} and
\ref{est4},   the fact that $u\to \delta_p$ as $t\to T^-$ and that
$h$ is smooth and bounded.
\end{proof}

\begin{cor}\label{LYHineq1} $v(x,t)\le 0$ on $M\times(0,T]$.
\end{cor}
\section{A pseudolocality theorem}
In this section, we will extend Perelman's pseudolocality theorem
to complete noncompact manifolds. We will prove the following:
 \begin{thm}\label{pseudolocality1} Let $n$ be fixed.
There exist $\delta, \e>0$ with the following property:

Suppose $g(x,t)$ is a smooth complete noncompact solution of the
Ricci flow with bounded curvature on $M^n\times [0,\e^2]$.
Suppose at some point $x_0\in M$ the isomperimetric constant in
$B_0(x_0,1)$ is larger than $(1-\delta)c_n$, where $c_n$ is the
isoperimetric constant of $\mathbb{R}^n$, and $R(x, 0)\geq -1$ for all $x\in B_0(x_0,1)$. Then $|Rm(x,t)|\le
    t^{-1}+\e^{-2}$ for $0<t\le \e^2$ and $x\in B_t(x_0,\e)$.
\end{thm}
By the result of \cite{Sh0}, we may   assume that the covariant
derivatives  of the curvature are uniformly bounded in spacetime.
The proof is similar to the case for compact manifolds using the
estimates obtained in previous sections. See \cite{P1,KL,CLN,ST}.
For the sake of completeness, we will sketch the proof.

Suppose this is not true. Then we can find $(M_i,g_i(t))$,
$\delta_i, \e_i>0$ with $\delta_i, \e_i\to0$ and $p_i\in M_i$
satisfying the following:
 \begin{enumerate}
    \item [\bf(b1)]$g_i(t)$ is a smooth solution of the Ricci flow on $[0,
     \e_i^2]$   with bounded $|\nabla^kRm|$ on
    $M_i\times[0, \e_i^2]$ for all $k\ge0$.
    \item [\bf(b2)] The isomperimetric constant in $B_0^{(i)}(p_i,1)$
is larger than $(1-\delta_i)c_n$.
    \item[\bf(b3)] There exist $0<t_i\le \e_i^2$, and $x_i\in B_{t_i}^{(i)}(p_i,\e_i)$
     and $|Rm(x_i,t_i)|\ge t_i^{-1}+\e^2_i$.
\end{enumerate}
Let $A_i=1/1000n\e_i$. By Claims 1 and 2 in \cite{P1}, see also
\cite{KL,CLN,ST}: We can find $\bar x_i,\bar t_i$ with $0<\bar
t_i\le \e_i^2$ and $\bar x_i\in B_{\bar t_i}^{(i)}(p_i,
(2A_i+1)\e_i)$ satisfying the following:

\begin{enumerate}
    \item [\bf(c1)]
$Q_i=|Rm(\bar x_i,\bar t_i)|\ge \frac1{\bar t_i}$, and if
$$
\bar t_i-\frac12Q_i^{-1}\le t\le \bar t_i,\ d_{\bar t_i}(x,\bar
x_i)\le \frac1{10}A_iQ_i^{-\frac12}
$$
then
$$
|Rm(x,t)|\le 4|Rm(\bar x_i,\bar t_i)|.
$$
\end{enumerate}

Consider the rescaled   flows: $\wh g_i(t)=Q_ig_i(\bar
t_i+Q_i^{-1}t_i).$ Then $\wh g_i$ satisfies the Ricci flow
equation  on $M_i\times[-\frac12,0]$ with bounded $|\nabla^kRm|$.
Moreover, the following are true:
\begin{enumerate}
    \item [\bf(d1)] $|Rm(\bar x_i,0)|=1$.
    \item [\bf(d2)] If $$
-\frac12\le t\le 0,\ d_{0}(x,\bar x_i)\le \frac1{10}A_i
$$
then
$$
|Rm(x,t)|\le 4.
$$
\end{enumerate}
Let $u_i$ be the fundamental solution of the conjugate heat
equation to the flow $\wh g_i$: $-(u_i)_t-\Delta u_i+R_iu_i=0$,
$\lim_{t\to 0}u=\delta_{\bar x_i}$. Let $v_i$ be the corresponding
LYH Harnack expression defined in (\ref{Harnackex}) with
$\tau=-t$. Then $v_i\le 0$ by Corollary \ref{LYHineq1}.

{\bf Case 1}: Suppose the injectivity radius at $\bar x_i$ at
$t=0$ are uniformly bounded from below. Because of {\bf (d2)} and
the Ricci flow equation, we know that the injectivity radius of
$\bar x_i$ are uniformly bounded from below at $t=-\frac12$. Since
$A_i\to\infty$, by the compactness result of  Ricci flow \cite{H1,
CLN}, for any sequence, we can find a subsequence of $\wh g_i$
which converge, still denoted by $\hat g_i$. Namely, there is
$(M,p,g(t))$ with $g(t)$ being a solution of the Ricci flow on
$[-\frac12,0]$, and an exhaustion  $U_i$ of $M$, diffeomorphisms
$\Phi_i:U_i\to M_i$ with the following properties:

\begin{enumerate}
    \item [\bf(e1)] $\Phi_{i}(p)=\bar x_i$.
    \item [\bf(e2)] $\Phi^*_{i}\hat g_i$ converges in $C^\infty$
    sense to $g$ on $M\times(-\frac12,0)$.
    \item[\bf(e3)] The curvature of $g(t)$ is bounded by 4.
    \item[\bf(e4)] There exists $0>t_0>-\frac12$ such that
    $|Rm(p,t)|\ge\frac12$ for all $t>t_0$.
\end{enumerate}

Note that {\bf (e4)} is a consequence of {\bf(d1), (d2), (e2)},
local derivatives bound for the curvature tensor and the evolution
equation of the curvature tensor.

\begin{lem}\label{converge1}\

\begin{enumerate}
    \item [\bf(f1)]   $\Phi_{i}^*u_i$
subconverge  on $M\times (-\frac12,0)$ to a solution of
$$
u_\tau-\Delta u+Ru=0.
$$
    \item[\bf(f2)] $u>0$, and if $v$  is the LYH Harnack
    expression defined in (\ref{Harnackex}) corresponding
    to $u$, then
$$
  v_\tau-\Delta v-Rv=-2\tau|R_{ij}+\nabla_i\nabla_j
  f-\frac1{2\tau}g_{ij}|^2
  $$
  where $f$ is given by $u=e^{-f}/(4\pi\tau)^\frac n2$.
 Moreover, $v\le0$.
\end{enumerate}
\end{lem}
\begin{proof}

 We first prove {\bf (f1)}.  By Lemma \ref{lm-int-bound}, $\int_{M_i}u_idV_t^{(i)}=1$ for all $i$ and
 $\tau$.
Since $u_i>0$ for $\tau>0$, by (\textbf{e2}) and the proofs of
Corollaries \ref{gradient-local} and
\ref{cor-mean-value-inequlity} , we conclude that $u_i$ are
locally uniformly bounded. By (\textbf{e2}), it is easy to see
that {\bf (f1)} is true. Note that we can construct the function
$\rho$ for $(M,g(t))$ as in the proof of Lemma \ref{lm-harnack}
and use this to prove a result similar to Corollary
\ref{gradient-local} for $u_i$  by (\textbf{e2}). See also Remark
\ref{gradientlocal1}.

Next we want to prove {\bf (f2)}. By the proof of Lemma
\ref{lm-int-lower-bound},  for a fixed but small neighborhood $U$
of $p$ there is $c>0$ and $\tau>0$, such that $\int_U
u_idV_{t}^{(i)}\ge c $ for all $i$.   So $u>0$ by the fact that
$\phi_i^*\wh g_i$ converges to $g$, {\bf (f1)} and the maximum
principle. The rest of the lemma follows from (\ref{Harnackex1}).
\end{proof}

 \begin{lem}\label{converge2} With the same notation as in Lemma \ref{converge1},
 for any
  $0<\tau_0< \frac12$ we have
  $$
  \int_{ B_{\tau_0 }(p,\sqrt{\tau_0})}vdV_{-\tau_0}<0.
  $$
 \end{lem}
\begin{proof} Suppose $\int_{ B_{\tau_0
}(p,\sqrt{\tau_0})}vdV_{-\tau_0}=0$, then $v=0$ in $B_{\tau_0
}(p,\sqrt{\tau_0})$. Let $h_0$ be a nonnegative smooth function
with support in $ B_{\tau_0 }(p,\sqrt{\tau_0})$ which is positive
somewhere.  Then for $i$ sufficiently large, we may also consider $h_0$ to be a
smooth function with compact support in $M_i$. Now solve the forward
heat equation with initial data $h_i|_{t=-\tau_0}=h_0$ on
$M_i\times [-\tau_0,0)$.  Then since $h_0$ is bounded, the $h_i$'s are also
uniformly bounded in space time. We may thus assume that $h_i\to h$ which solves the heat equation in $(M,g(t))$ and the convergence is
uniform on compact sets of $(x,t)\in M\times[-\tau_0, 0)$. By
Theorem \ref{LYHineq}(ii), for $0<\tau<\tau_0$

$$
\int_{M_i}v_i(x,-\tau)h_i(x,-\tau)dV_{-\tau}^{(i)}\ge
\int_{M_i}v_i(x,-\tau_0)h_0(x)dV_{-\tau_0}^{(i)}
$$
where $v_i$ is the LYH Harnack expression for $v_i$. Since
$v_i\le0$ by Corollary \ref{LYHineq1} and since $h_0$ is a fixed
function with compact support, let $i\to\infty$, we can  conclude
for any compact set $K$ in $M$,
$$
\int_{K}v(x,-\tau)h(x,-\tau)dV_{-\tau} \ge \int_{M }v
(x,-\tau_0)h_0(x)dV_{-\tau_0}=0
$$
 Since $v\le0$ and $h>0$ for $t>-\tau_0$, we have $v=0$ for
$0<\tau<\tau_0$. By {\bf (f2)} we have
$$
R_{ij}+f_{ij}-\frac1{2\tau}g_{ij}=0.
$$
 for $\tau<\tau_0$. Since the curvature is uniformly bounded, for any $0<\tau<\tau_0$,
 $|\nabla f|$ is at most linear growth. From this one can prove
   that the
 vector field  $Y_t=(1-\frac1{ \tau }t)^{-1}\nabla f(-\tau )$ can be
 integrated from $0$ to $t$ as long as $1-\frac1{ \tau}t>0$
 which defines a diffeomorphism $\psi_t$.

 Then the flow $\wt g(t)=(1-\frac1{ \tau_1} t)\psi_t^*(g(-\tau_1))$
 is a solution of the Ricci flow on $[0,\tau_1)$ with initial data
 $g(-\tau_1)$ \cite[p. 22-23]{CK}. Now $g(t-\tau)$ for
 $0\le t<\tau$ is also such solution. Since the curvatures are
 bounded for both flows, by the uniqueness result of \cite{CZ},
 they are the same. However, by {\bf (e4)} the curvature of $(1-\frac1{
 \tau}
 t)\psi_t^*(g(-\tau_1))$ blows up near $t=\tau$ and the curvature
 of $g(t-\tau)$ are uniformly bounded, so this is impossible.
 \end{proof}

{\bf Case 2}: Suppose the injectivity radii $a_i$ at $\bar x_i$ at
$t=0$ tend to zero. We further rescale the metrics: Let $ \Hat
 {\Hat g}_i(t)=a_i^{-1}\hat g_i(a_i t)$. Then it is defined on
$[-\frac12a_i^{-1},0]$ with the following properties:

\begin{enumerate}
    \item [\bf(g1)]  $|Rm(\bar x_i,0)|=a_i$.
    \item [\bf(g2)] If $$
-\frac12a_i^{-1}\le t\le 0,\ d_{0}(x,\bar x_i)\le
\frac1{10}A_i\cdot a_i^{-\frac12}
$$
then
$$
|Rm(x,t)|\le 4a_i.
$$
\item [\bf(g3)] There exists $t_0<0$ such that the injectivity
radius of $ \bar x_i$ at time $t_0$ is less than 2.
\end{enumerate}

{\bf(g3)} can be proved by the fact that $a_i\to0$ and the
injectivity radius bound in \cite{CGT}.

 As before,
we can find a limit metric and flow $(M,g(t),p)$ on $(-\infty,0]$
with the following properties:

\begin{enumerate}
    \item [\bf(h1)] $g(t)$ is flat for all $t$.
    \item [\bf(h2)] The injectivity radius of $p$
    at time $t_0$ is less than or equal to $2$.
\end{enumerate}

Let $v$ and $f$ as before.

\begin{lem}\label{converge4} For any
  $0<\tau_0< \infty$ we have
  $$
  \int_{ B_{\tau_0 }(p,\sqrt{\tau_0})}vdV_{-\tau_0}<0.
  $$
 \end{lem}
\begin{proof} As before, if this is not true, then one can prove
that
$$
R_{ij}+f_{ij}-\frac1{2\tau}g_{ij}=0.
$$
for $0<\tau<\tau_0$. Since curvature is bounded, there is
$0<\tau_1<\tau_0$ such that $f_{ij}\ge g_{ij}$. One can prove that
\begin{equation}\label{converge5}
 g(t-\tau_1)=(1-\frac1{ \tau_1} t)\psi_t^*(g(-\tau_1)).
\end{equation}

The function $f$ is an exhaustion function, see \cite{CT3} for
example. So there is a point such that $\nabla f=0$. Hence we can
find a point $q$ which is a fixed point of $\psi_t$. Hence the
injectivity radius of $q$ with respect to $\psi_t^*(g(-\tau_1))$
is independent of $t$. Note that $M$ is flat but is not
$\mathbb{R}^n$ by {\bf (h2)}, and so the
  injectivity radius of $q$ is finite. On the
other hand, $g(t)=g(-\tau_1)$ because of {\bf(h1)} and the fact
that $g$ satisfies the Ricci flow equation. Hence the injectivity
radius of $q$ of $g(t)$ is independent of $t$. This is impossible
as $t\to \tau_1$ by (\ref{converge5}).
\end{proof}
From these, it is easy to see that Claim 3 in \cite[\S10]{P1} is
also true. One can conclude that Theorem \ref{pseudolocality1} is
true by the argument of \cite{P1}, see also \cite{KL,CLN,N1,ST}.

 \section{singularity formation and longtime existence}

We now apply the pseudolocality result to describe where
singularities to the Ricci flow can form under certain
assumptions.  More precisely, in Theorem \ref{ltt1} we prove that
any finite time singularities of the Ricci flow (\ref{lte1}) under
Assumption \ref{ltl1} below must form within a compact
set. In Theorem \ref{ltt2} we apply this result to complete
non-negatively curved \K manifolds and prove a long time existence
result for the \KR flow.

Consider the Ricci flow
\begin{equation}\label{lte1}
  \begin{cases}\frac{d}{dt}g&=-2Rc\\
  g(t)&=g.
  \end{cases}
  \end{equation}

on a complete non-compact Riemannian manifold.  Let us make the
following assumption

\begin{ass}\label{ltl1} Let $(M, g)$ be a complete non-compact
Riemannian manifold of dimension $n$ such that
\begin{enumerate}
 \item $|Rm(x)|\to 0$ as
$d(x, p)\to\infty$ for some fixed point $p$ \item the injectivity
radius $inj(M, g)$ of $(M, g)$ is bounded from below.
\end{enumerate}
\end{ass}

\begin{lem}\label{ltl2}
Let $(M, g)$ be as in Assumption \ref{ltl1} and let $(M, g(t))$ be
the corresponding $maximal$ solution to the Ricci flow
(\ref{lte1}) on $M\times[0, T)$.  For any $0< \delta < 1$, there
exists $r
>0$ with the following property: Given any $t' <T$ there exists
$0<d'<\infty$ such that for any $x\in M$ and $0<t\le t'$  with
$d_{t}(x, p) \geq d'$, we have $$Vol_{t}(\partial \Omega)^n \geq
(1-\delta) c_n Vol_{t}(\Omega)^{n-1}$$ for any $\Omega \subset
B_{t}(x, r)$.
\end{lem}
\begin{proof}\hspace{12pt}
Let $0< \delta < 1$ be given and let $r_0= {\text{inj}(M,
g(0))}/{2}$. By conditions (1) and (2) in Assumption \ref{ltl1},
we can find some $0<d_0 <\infty$ such that for any $x$ where
$d_{0}(x, p) \geq d_0$, we have
\begin{equation}\label{lte2}
Vol_0(\partial \Omega)^n \geq (1-\frac{\delta}{2}) c_n
Vol_0(\Omega)^{n-1}
\end{equation}
for any $\Omega \subset B_{0}(x, r_0)$.

By Theorem 18.2 in \cite{H2}, given   $t' <T$ and   any $\eta
>0$ there exists some compact $S \subset M$ such that
\begin{equation}\label{lte3}
|Rm(x, t)| \leq \eta
\end{equation}
for all $(x, t) \in (M \setminus S) \times [0, t']$.

By choosing $\eta$ sufficiently small we see from (\ref{lte1}),
(\ref{lte2}) and (\ref{lte3}) that we may choose some $d'>d_0$
such that for any $x$ where $d_{t}(x, p) \geq d'$, we have
\begin{equation}\label{lte4}
Vol_{t}(\partial \Omega)^n \geq (1-\delta) c_n
Vol_{t}(\Omega)^{n-1}
\end{equation}
for any $\Omega \subset B_{t}(x, \frac{r_0}{2})$.  Thus
$\frac{r_0}{2}$ satisfies the conclusion of the Lemma.
\end{proof}

\begin{thm}\label{ltt1}
Let $(M, g)$ satisfy Assumption \ref{ltl1} and let $(M, g(t))$ be
the corresponding $maximal$ solution to the Ricci flow
(\ref{lte1}) on $M\times[0, T)$.  Then either $T=\infty$ or there
exists some compact set $S \subset M$ with the property that
$|Rm(x, t)|$ is uniformly bounded on $(M \setminus S) \times [0,
T)$.
\end{thm}

\begin{proof}
Assume that $T< \infty$ and let $g(t)$ be a $maximal$ solution to
(\ref{lte1}) on $M\times[0, T)$.  Thus there exists sequences $x_i
\in M$ and $t_i \to T$ such that $|Rm(x_i, t_i)| \to \infty$ as
$i\to \infty$.  We will show that there exists some compact $S
\subset M$ such that every such sequence $x_i$  must be contained
inside $S$. $S$ will then clearly satisfy the conclusion of the
Theorem.

 Suppose there is a sequence $(x_i,t_i)$ satisfying the above condition and
   $d_0(p, x_i) \to \infty$.  Let $\delta, \epsilon$
be as in Theorem \ref{pseudolocality1}.
 For such  $\delta$, let $r$ be as in Lemma
\ref{ltl2}.  By rescaling
 our solution $g(t)$ in both time and space, we may assume that $r=1$
  (without affecting $\delta$).
Now let $t'=T- \epsilon^2 $, and choose $d'$ as in Lemma
\ref{ltl2}. Then by Theorem 18.2 in \cite{H2}, we may assume $d'$ sufficiently large so that $R(y, t') \geq =-1$ for all $y\in B_{t'}(y, x)$ where $d_{t'}(x, p)\geq 0$. We may assume that $\e>0$ is small enough such
that $t'>0$.

Let $\eta_k\to0$ and let $\tau_k=t'-\eta_k>0$ and
$g_k(t)=g(\tau_k+t)$. Then $g_k(t)$ is well defined on $[0,\e^2]$.
By \cite{Sh0,H2}, we know that for each $k$ the curvature tensor of $g_k(t)$
together with its derivatives are uniformly bounded in $[0,\e^2]$.
Let $i_0$ be large enough such that if $i\ge i_0$, then
$d_{t}(x_i,p)\ge d'$ for all $0<t\le t'$.

 By Theorem \ref{pseudolocality1}, we have
$$
|Rm_k(x_i,t)|\le t^{-1}+\e^{-2}
$$
for all $i\ge i_0$ and $0\le t\le \e^2$. Here $Rm_k$ is the
curvature tensor for $g_k$. Now
$Rm(x_i,t_i)=Rm_k(x_i,t_i-\tau_k)$.  We have
$$
t_i-\tau_k=t_i-t'+\eta_k=t_i-T+\e^2+\eta_k.
$$
Hence for fixed $i\ge i_0$ such that $T-t_i \le \e^2$, we have $0\le t_i-\tau_k\le \e^2$ for $k$ sufficiently large.

So
$$
|Rm(x_i,t_i)|=|Rm_k(x_i,t_i-\tau_k)|\le (t_i-\tau_k)^{-1}+\e^{-2}.
$$
Now letting $k\to\infty$ and then letting $i\to\infty$, we have
$$
\limsup_{i\to\infty}|Rm(x_i,t_i)|\le 2\e^{-2}.
$$
This contradicts our initial assumption, and thus completes the proof of the Theorem.

\end{proof}

\begin{cor}
Suppose $T < \infty$ in Theorem \ref{ltt1}.  Then $Rm(x, T) \to 0$ as $x\to \infty$ in the sense that: given any $\epsilon>0$, we may choose $S$ such that $|Rm(x, t)| \leq \epsilon$ for all $(x, t)\in S^c \times [0, T)$.
\end{cor}
\begin{proof}

Assume the Corollary is false.  Thus there exits a space time sequence $(x_k, t_k)$ such that $d_0(p, x_k) \to \infty$, $t_k \to T$ and $|Rm(x_k, t_k)|\geq C_1$ for some $C_1>0$.

\vspace{12pt}

1. Fix some small $\epsilon _1 >0$ to be chosen later and let $s_k=t_k-\epsilon_1$.  Then by Theorem \ref{ltt1} and Theorem 18.2 and 13.1 in \cite{H2}, we may assume the compact set $S \subset M$ from Theorem \ref{ltt1} was chosen sufficiently large so that for some $C_2 >0$ we have $$|\frac{d}{dt}Rm(x_k, t)|\leq C_2$$ for all $k$ sufficiently large and $t \in [0, t_k]$.
 \vspace{12pt}

2. For $k$ suffuciently large, $s_k \in [0, T-\epsilon_1]$.  Thus by Theorem 18.2 in \cite{H2}, given any $\epsilon_2 >0$ we may assume $S \subset M$ was chosen sufficiently larger still so that $$|Rm(x_k, s_k)|\leq \epsilon_2 $$ for $k$ sufficiently large.
 \vspace{12pt}

Thus by 1 and 2, for $k$ sufficiently large we have that 
$$|Rm(x_k, t_k)|\leq |Rm(x_k, s_k)|+C_2 \epsilon_2 \leq \epsilon_2+C_2 \epsilon_1 .$$

Note that $\epsilon_1$, $\epsilon_2$ were chosen independently of each other, and that $C_2$ can be chosen independent of these.  Thus for sufficiently small choices of $\epsilon_1$ and $\epsilon_2$, we arrive at a contradiction.  This compeletes the proof by contradiction.
\end{proof}

We now apply Theorem \ref{ltt1} to the case non-negatively curved
\K \ manifolds.

\begin{thm}\label{ltt2}
Let $(M, g)$ be a complete non-compact \K \ manifold with
non-negative holomorphic bisectional curvature, strictly positive
at some point.  Then if $(M, g)$ satisfies Assumption \ref{ltl1},
the \KR flow has a long time solution $g(t)$ on $M\times [0,
\infty)$.
\end{thm}
\begin{proof}
Assume the Theorem is false and that $g(t)$ is a $maximal$
solution on $M\times[0, T)$ for $T< \infty$.  Let the set $S$ be
as in Theorem \ref{ltt1}.  Now define $$F(x, t)=\log \frac{\det
g(x,t)}{\det g(x,0)}.$$ By \cite{Sh2}, we know that $g(t)$ also
has non-negative holomorphic bisectional curvature and hence
$F(x,t)\le0$. We claim that $-F(x,t)$ is uniformly bounded on
$M\times[0,T)$. If this is true, then the curvature is also
uniformly bounded on $M\times[0,T)$ by the argument in
\cite[\S7]{Sh2}. From this it is easy to see that the solution
$g(t)$ can be extended beyond $T$ by \cite{Sh0,CZ}. This is a
contradiction.

Since the curvature  $|Rm(x,t)|$ is uniformly bounded in
$(M\setminus S)\times [0, T)$ and since
$$
-F(x,t)=\int_0^t R(x,\tau)d\tau
$$
where $R$ is the scalar curvature,
 there exists $C_1$ such that
\begin{equation}\label{Fbound}
0\le - F(x, t) \leq C_1
\end{equation}
for all $(x, t) \in (M\setminus S)\times [0, T)$.

Next we want to prove that  there exists $C_2$ such that
\begin{equation}\label{Fbound1}
-F(x, t)\leq C_2
\end{equation} for all $(x, t) \in S\times [0, T)$.

Let $\tilde M$ be the universal cover of $M$. Then by \cite{NT1},
$\tilde M=\tilde N\times \tilde L$ holomorphically and
isometrically where $\tilde N$ is compact and $\tilde L$ satisfies

\begin{equation}\label{scalarcurvaturedeca}
\frac{1}{\tilde V_{\tilde o}(r)}\int_{\tilde B_{\tilde
o}(r)}\tilde R\le \frac{C}{1+\tilde r}
\end{equation}
where $\tilde B_{\tilde o}(r)$ is the geodesic ball in $\tilde L$
and $\tilde V_{\tilde o}(r)$ is its volume. Also $\tilde R$ is the
scalar curvature of $\tilde L$. Since $|Rm(x)|\to 0$ as $x\to
\infty$ in $M$, we must have $\tilde M=\tilde L$.

 Suppose  there exist  sequences $x_i \in S$ and $t_i\to T$ such that
$F(x_i, t_i)\to -\infty$. By (\ref{Fbound}), we may assume that
$$m(t_i):=\min_{M}F(\cdot, t_i)=F(x_i, t_i)$$ provided $i$ is
sufficiently large. After we lift the K\"ahler-Ricci flwo to
$\tilde M$,   there is a compact set $\tilde S$ in $\tilde M$ and
$\tilde x_i\in \tilde S$ such that $F(\tilde x_i, t_i)\to -\infty$
and $\min_{\tilde M}F(\cdot, t_i)=F(\tilde x_i, t_i)=m(t_i)$. Here
$F(\tilde x,t)$ is the logarithm of the ratio of volume elements
for the flow in $\tilde M$.

Then by the proof of Corollary 2.1 of \cite{NT2}  we  have
\begin{equation}\label{lte5}
\begin{split}
-m(t_i)\leq & C_3 \int_{0}^{\sqrt{at_i (1-m(t_i))}} \frac{s}{1+s}   ds\\
  \leq & C_3  \sqrt{at_i (1-m(t_i))}
  \end{split}
\end{equation}
for some positive constants $a, C_3$, where we have used the fact
that $\tilde S$ is compact and (\ref{scalarcurvaturedeca}). From
this we can see that (\ref{Fbound1}) is true. This completes the
proof of the claim and the theorem.
\end{proof}

\end{document}